\documentclass{aims}
\usepackage{amsmath}
  \usepackage{paralist}
  \usepackage{graphics} %% add this and next lines if pictures should be in esp format
  \usepackage{epsfig} %For pictures: screened artwork should be set up with an 85 or 100 line screen
\usepackage{graphicx}  \usepackage{epstopdf}%This is to transfer .eps figure to .pdf figure; please compile your paper using PDFLeTex or PDFTeXify.
 \usepackage[colorlinks=true]{hyperref}
   \usepackage{amsthm,amscd,amsxtra,amsfonts,amsmath,amssymb,multirow,mathrsfs,
url}
\DeclareMathOperator*{\argmin}{arg\,min}
\usepackage{algorithm}
\usepackage{algpseudocode}
\usepackage{float}
\usepackage{multicol}
\usepackage{listings}
\usepackage{verbatim}
\usepackage{wrapfig,epsfig}
\usepackage{amsmath}
\usepackage{amsfonts}
\usepackage{amssymb}
\usepackage{graphicx}%
\usepackage{titlesec}
\usepackage{dsfont}
\usepackage{physics}
\hypersetup{urlcolor=blue, citecolor=red}

\def\currentvolume{X}
 \def\currentissue{X}
  \def\currentyear{200X}
   \def\currentmonth{XX}
    \def\ppages{X--XX}

\newtheorem{theorem}{Theorem}[section]

\theoremstyle{definition}
\newtheorem{definition}[theorem]{Definition}

\newcommand{\R}{\mathbb R}

\newcommand{\bd}[1]{\vb*{#1}}

\title[Geometric structures guided model] %Use the shortened version of the full title
      {Geometric structure guided model and algorithms for complete deconvolution of gene expression data }

\author[Duan  Chen, Shaoyu Li and Xue Wang]{}


\thanks{$^*$ Corresponding author: Duan Chen}

\begin{document}
\maketitle

% Enter the first author's name and address:
\centerline{\scshape Duan Chen$^1*$}
\medskip
{\footnotesize
% please put the address of the first author
 \centerline{Department of Mathematics and Statistics}
   %\centerline{Other lines}
   \centerline{ University of North Carolina at Charlotte, USA}
} % Do not forget to end the {\footnotesize by the sign }

\medskip

\centerline{\scshape Shaoyu Li$^1$ and Xue Wang$^2$}
\medskip
{\footnotesize
 % please put the address of the second  and third author
 \centerline{$^1$Department of Mathematics and Statistics}
  \centerline{University of North Carolina at Charlotte, USA}
%   \centerline{Other lines}
   \centerline{$^2$Department of Quantitative Health Sciences, }
    \centerline{Mayo Clinic, Florida, 32224}
%   \centerline{Other lines}
  
}

\bigskip

% The name of the associate editor will be entered by an editorial staff
% "Communicated by the associate editor name" is not needed for special issue.
 %\centerline{(Communicated by the associate editor name)}

%The abstract of your paper
\begin{abstract}
Complete deconvolution analysis for bulk RNAseq data is important and helpful to distinguish whether the difference of disease-associated GEPs (gene expression  profiles) in tissues of patients and normal controls are due to changes in cellular composition of tissue samples, or due to GEPs changes in specific cells.
One of the major techniques to perform complete deconvolution is nonnegative matrix factorization (NMF), which also has a wide-range of applications  in the machine learning community. However, the NMF is a well-known strongly ill-posed problem, so a direct application of NMF to RNAseq data will suffer severe difficulties in the interpretability of solutions.  In this paper we  develop  an NMF-based mathematical model and corresponding  computational  algorithms to improve the solution identifiability of deconvoluting bulk RNAseq data. In our approach, we combine the biological concept of marker genes  with the solvability conditions of the NMF theories, and  develop a geometric structured guided optimization model. In this strategy,  the geometric structure of bulk tissue data is first explored by the spectral clustering technique. Then,  the identified information of marker genes  is integrated as solvability constraints, while the overall correlation graph is used as manifold regularization. Both synthetic and biological data are used to validate the proposed model and algorithms, from which solution interpretability and accuracy are significantly improved.

\end{abstract}

\section{Introduction}

	Over past decades, analysis of   transcriptome or gene expression data has been an essential component to understand the processes involved in human development and disease \cite{cang2020inferring,jin2020scai,zhang2019revealing,harrington2017geometric}, but  the complex nature of tissue samples   under investigation remains as a major obstacle\cite{davey1996flow,whitney2003individuality,de2005purity}. A bulk tissue sample could include many cell types, and its heterogeneous characteristics make the interpretation of gene expression (such as RNA-Seq) complicated \cite{fridman2012immune,shen2013computational,avila2018computational}: 
	for every gene, its  measured gene expression profiles (GEPs)  in a compound  sample are actually tissue-averaged, i.e.,  the sum of expression of   all cells in the sample. On the other side,  cellular composition of bulk samples varies, and different samples may show high variance between one and another in relative cell subset proportions. So the GEP of low abundant cell types could be masked by that of ones with higher proportions. Consequently, 
it is challenging to  determine whether an experimental or clinic treatment should target one particular {gene type}  or focus on investigating possible sources of varying { cell types}  among samples. 
 For example, Alzheimer's disease (AD) is marked 
	 by amyloid-beta plaques and neurofibrillary tangles, along with neuronal loss and gliosis in the affected brain regions. 
	Transcriptome-wide GEP from brain tissue  of AD patients and neuropathologically normal controls are different. Such differences are  critical for discovering genes and biological pathways that are perturbed in and/or lead to AD \cite{allen2018conserved,mckenzie2017multiscale,mostafavi2018molecular,de2018multi}. 
	 Differential expression (DE) analysis is one of the important tools to unveil these differences. It will reveal novel insights into the genes and pathways, and is potentially helpful for drug targets  AD therapeutics.
	  However, a fundamental knowledge gap still remains for DE, concerning  whether disease-associated GEP changes in brain   tissues are due to changes in cellular {\em composition} of tissue samples, or due to GEP changes in {\em specific cells}, e.g., central nervous system (CNS) cells.
	   It could be much more informative to  study gene expression on specific cells, or  identify cell-intrinsic differentially expressed genes (CI-DEGs).
	     But for many complex biological mixtures, exhaustive knowledge of individual cell types in brain tissues and their specific markers is lacking. 
	     Although single-cell RNA sequencing (sRNAseq) data can be used or   serve as a reference, such approaches remain costly, cumbersome and limited in sample sizes\cite{zhang2016purification,darmanis2015survey,lake2018integrative}. 
  \begin{figure}[h!]
		\begin{center}
			\includegraphics[width=0.65\textwidth]{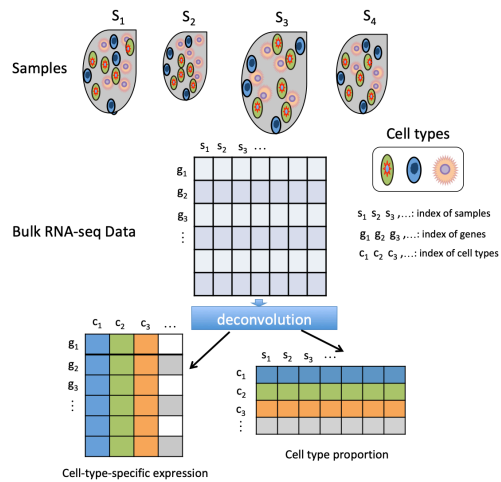} 
		\end{center}
		\caption{Diagram of complete deconvolution of bulk tissue data}%
		\label{fig:decov}%
	\end{figure}

	     In contrast, computational tools can be used to leverage widely available large-scale bulk tissue RNAseq data sets \cite{mckenzie2017multiscale,de2018multi,allen2016human,kuhn2011population,chikina2015cellcode}. This problem, as illustrated in Figure \ref{fig:decov}, is  called complete deconvolution. In this approach, expression of a gene in a sample tissue is assumed to be the linear combination of its expressions in the constituting cell types, with respected to the cell proportion, i.e.
	     \begin{equation}\label{eqn:basic}
	     	g_{ij} = \sum_{l=1}^kp_{lj}c_{il}, \quad 1\le i \le N, \quad 1\le j\le n,
	     \end{equation}
	     where $g_{ij}$ and $c_{i,l}$ are the GEPs of gene $i$ in the $j$-th sample and $l$-th cell type, respectively, while $p_{lj}$ is the proportion of the $l$-th cell type in the $j$-th sample. For the total number of genes $N$, total number of samples $n$, and the number of cell types $k$, we usually have $N\gg n>k$. In matrix form, Eq. (\ref{eqn:basic}) is represented as ${\bf G} = {\bf CP}$ with all matrix entries being non-negative. Given data ${\bf G}$, {\em both} variables $\bf C$ and $\bf P$ are to be solved.  Note that many  deconvolution  algorithms have been developed \cite{tsoucas2019accurate,avila2018computational,mohammadi2016critical,newman2015robust,qiao2012pert,zhong2013digital,gong2013deconrnaseq,cui2016gene,abbas2009deconvolution, gaujoux2012semi, shen2013computational} for GEP  in bulk tissues, but their primary focuses have been only on estimating the cellular composition with prior known information of cell-specific markers. This type of problem to solve for $\bf P$, with known $\bf C$, is called partial convolution, can be performed with remarkable robustness and accuracy. However, in more realistic circumstances when little or no information about the underlying cell type is available, developing reliable complete deconvolution methods is still an open problem and only a handful models have been established \cite{kang2019cdseq,repsilber2010biomarker,zaitsev2019complete}.   
	     %On the other hand, estimating and benchmarking GEP in specific cell types is quite important, such as in co-expression network analysis and expression quantitative trait loci (eQTL) analysis. 
	     
	     Mathematically, complete deconvolution  is a  nonnegative matrix factorization (NMF) problem \cite{craig1994minimum,paatero1994positive,lee1999learning}. Many studies  have been  established  for various types of data in  other fields, such as spectral unmixing in analytical chemistry \cite{paatero1994positive}, remote sensing \cite{ma2013signal}, image processing \cite{fu2015self}, or topic mining in machine learning \cite{zhang2006learning}, etc.
    Fundamental computational algorithms include the Multiplicative Update Algorithms (MUA)\cite{lee1999learning} and alternating nonnegativity constrained least squares (ANLS)\cite{kim2008nonnegative}.
There is no obstacle at all if one is simply looking for a couple of solutions $\bf C$ and $\bf P$. However, the NMF is strongly ill-posed and  solutions  are generally not separable, i.e., for any ${\bf \Omega}\in\R^{k\times k}$, $\hat{\bf C}={\bf C}{\bf \Omega}$ and  $\hat{\bf P}={\bf \Omega}^{-1}{\bf P}$ are also solutions, as long as their non-negativity is satisfied. Such non-uniqueness will  pose great challenges on  {\em solution interpretability}: For RNAseq data,  different solutions represent various combinations of GEPs in each cell types and cell proportions in tissues.  Meaningful explanation of these biological quantities is critical to next step DE analysis. 
There are a few guidelines to reduce such ill-posedness. As stated in, if the matrices $\bf C$ and $\bf P$ satisfies certain identifiability conditions or  structures, see the details in Section \ref{sec:nnmf}, it is possible to have unique solutions,  subjective to row/column scaling and permutation ambiguities. 
Applying these sufficient conditions depends on the specific properties of available data in the corresponding research field. 
%Techniques addressing this ill-posedness heavily depend on the specific data characteristics in various applications.  
Successful methods in one field cannot be directly implanted to another because of different data characteristics. Modeling the right NMF tool for the application at hand is essential \cite{fu2019nonnegative}.  
On the biology side, the GEP data  $\bf G$  in bulk tissues includes expression of  marker genes, or cell-type-specific genes, which are defined by their exclusive expression in only one component (cell type) in cell mixtures. In the ideal, noise-free scenario, this property implies that GEP of  a marker gene of across samples will be exactly linearly dependent to the proportion of the corresponding cell type across the samples.  In the realistic case, GEP vectors across tissue samples for all marker genes of the {\em same } cell type will display strong correlations, while strong orthogonality for  {\em different} cell types.  These biological characteristics establish a connection to the mathematical theories of NMF. So it is possible to develop robust and accurate complete deconvolution algorithms for bulk tissue GEP data, without prior information of GEP in single cells.
%But before using NMF as a model to perform RNAseq data analysis, one needs to answer the following questions:  (a) how to recognize the structure details only from the available data? (b) how to use these sufficient conditions as optimization constraints to the NMF problem?

	  The   objective of the current work is to develop mathematical model and computational algorithms to perform complete deconvolution of bulk RNAseq data with reduced solution ambiguity and hence high interpretability. Our approaches are based on the abovementioned inherent characteristics of bulk tissue data and the theoretical foundation of NMF problem.
	%A common feature of these bulk tissue   data is that many GEP across tissue samples, or rows of $\bf G$,  exhibit strong mutual correlations.  This is due to the exclusive expression of cell type-specific genes, or marker genes,  in only one component within a mixture. This characteristics shapes a special geometric structure of the biological data, thus it is promising for the unique solution of the NMF problem. 
	And the goal is achieved by a structure-exploring and inheriting strategy. To explore the structure, we first define correlation distance among rows of data $\bf G$ (considered as data points in $\R^n$) and generate the corresponding graph, then spectral clustering technique is used to classify all points in $k$ (assumed number of cell types) groups. Finally, marker genes of each cell type are identified by picking the most correlated points in each cluster. Note that, in the noiseless case, rows of $\bf C$ can be understood as coefficients of data points with rows of $\bf P$ as coordinates.  Thus, we expect the row space of $\bf C$ inherits the geometric structure of $\bf G$ and impose the weak identifiability (to accommodate noises in real data) condition on $\bf C$.  Additionally, the manifold regularization is applied by the local invariance assumption \cite{belkin2002laplacian,cai2009probabilistic,he2004locality}. Combining these approaches, we establish a structure guided non-convex optimization model, to deconvolute bulk tissue RNAseq data without a prior information about marker genes. The proposed model is numerically solved under the frame work of alternating direction method of multipliers (ADMM)\cite{eckstein2012augmented,boyd2011distributed}, in which each variable can be solved one at a time in a two-fold iteration. Effectiveness and accuracy of the model and algorithms are tested by both synthetic and biological data. 
	This work is motivated by various manifold regularization NMF models, such as\cite{cai2010graph,qin2019fast}, but it has the following novel features: (1) Traditionally, Euclidean  distance is used in the manifold assumption of data space and it results a linear graph regularization term. But the application on biological data requires  the  correlation distance, from which a nonlinear graph regularizor is derived and it poses great challenges in computation; (2) More importantly, this new model is equipped with a solvability constraint, and this regularizor significantly improves solution identifiability from realistic noisy data.

%	  : use biological characteristics to enforce solvable conditions to enhance interpretability. 1, finding the structure by clustering; 2. preserve the structure by constraints; 3. computational algorithms; 4. numerical experiments; 5. possible challenges.
	     
	     The paper is organized as the following: Section \ref{sec:nnmf} briefly reviews the NMF and its separability conditions, and why this condition is related to the biological problem. 
	     Section \ref{sec:model} presents the geometric structure guided complete deconvolution model, including using spectral clustering analysis to identify marker genes (finding structures) and the quantitative constraints in the optimization problem (preserving the structure). For the resulting non-convex learning model, an  ADMM based algorithm is introduced in Section \ref{sec:algorithm}. As validations, in Section \ref{sec:num} there display numerical results of the proposed model and algorithms for both synthetic and biological data. The paper ends with a conclusion in Section \ref{sec:con}, where potential challenges of the work and possible future research directions are discussed.

\section{NMF and its identifiability conditions} \label{sec:nnmf}
In this section, we briefly review notations and some theoretical foundations of  NMF.

\subsection{Notations} Throughout the paper, a bold lower-case letter, such as $\bd{x}$, represents a column vector with the appropriate dimension, and $|\bd{x}|$ represents its $l_2$ norm. Vector $\bf 1$ is a column vector of some dimension with all entries being one. For a  matrix $\bf A$, $\|{\bf A}\|_F$ represents  its Frobenius norm,  ${\bf A}_{(i)}$ and ${\bf A}^{(j)}$ mean its $i$-th row and $j$-th column, respectively. 

Let ${\bf G}\in\R^{N\times n}$ with entry $g_{ij}$ be the expression of the $i$-th gene in the $j$-th sample; ${\bf C}\in\R^{N\times k}$ with entry $c_{ij}$ being the reference expression of the $i$-th gene in the $j$-th cell type; and ${\bf P}\in\R^{k\times n}$ with entry $p_{ij}$ being the proportion of the $i$-th cell type in the $j$-th sample. Dimensions $N \gg \max{(n,k)}$. The following linear relation is assumed: 
\begin{equation}
	{\bf G}={\bf CP} + \epsilon,
\end{equation}
where $\epsilon$ is noise. The  problem of complete deconvolution can be summarized as:  given data ${\bf G}\in\R^{N\times n}$, solve 
\begin{equation}\label{eqn:dconv}
	({\bf C}^*, {\bf P}^*) = \argmin_{{\bf C}\in\R_+^{N\times k}, {\bf P}\in\R_+^{k\times n}}{\delta({\bf CP, G})} 
\end{equation}
where $\R_+^{N\times k}$ or $\R_+^{k\times n}$ represent matrices with nonnegative entries and $\delta(\cdot, \cdot)$ is a cost function. For simplicity, we consider $\delta({\bf CP, G})=\displaystyle{\frac{1}{2}\|{\bf G - CP}\|_F^2}$ in this paper. 

\subsection{Ill-posedness of NMF:} Solving Eq. (\ref{eqn:dconv}) for only $\bf C$ (or $\bf P$) with the other variable known (partial deconvolution) is simply a convex regression problem. But it is well-known that solving both variables simultaneously is non-convex, NP-hard in general, and computational algorithms only converge to local minima or just stationary points \cite{wang2012nonnegative}.  Further, the NMF is ill-posed and the solution is not unique, or identifiable: if $({\bf C}^*, {\bf P}^*)$ is a local minimum to (\ref{eqn:dconv}), then for any ${\bf \Omega}\in\R^{k\times k}$, $\hat{\bf C}={\bf C}^*{\bf \Omega}$ and  $\hat{\bf P}={\bf \Omega}^{-1}{\bf P}^*$ are also solutions, as long as their non-negativity is satisfied. Non-uniqueness of solution will significantly impact statistical analysis for decisions in biological implementation. Therefore, it is important to restrict searching space of variables to increase the identifiability of solutions, in order for better interpretability. First,  uniqueness   of NMF solution is defined in the following sense \cite{fu2019nonnegative}:
\begin{definition} [Uniqueness (indentifiability) of NMF solution] 
The solution $({\bf C}^*, {\bf P}^*)$ of NMF (\ref{eqn:dconv}) is unique, or identifiable, if and only if for any other solution $(\bar{\bf C}, \bar{\bf P})$, there exists a permutation matrix $\Pi\in\{0,1\}^{k\times k}$ and a diagonal scaling matrix $\bf S$ with positive diagonal matrix such that
	\begin{equation}\label{eqn:def}
		\bar{\bf C} = {\bf C}^*\Pi{\bf S} \quad \text{ and } \quad \bar{\bf P} = {\bf S}^{-1}\Pi^\top{\bf P}^*. 
	\end{equation}
\end{definition}

\subsection{Geometric interpretation.}
It is summarized in \cite{wang2012nonnegative,fu2019nonnegative} that the uniqueness can be achieved under certain circumstances. In order to understand the strong and weak conditions of uniqueness, we  review the NMF problem from the perspective of geometric structures.
For  ${\bf A}\in\R^{m\times n}$,  the notation $cone({\bf A})$ denotes  the convex cone generated by the columns of ${\bf A}$, i.e.
\begin{equation}
	{cone}\{\mathbf{A}\}=\{\bd{x}\in\R^m|\bd{x}=\mathbf{A}\theta, \quad\text{ for some }\theta\in\R^n,  \theta\ge \mathbf{0}\},
\end{equation}
 and the conex hull of $\bf A$ is defined as 
 \begin{equation}
	 {conv}\{{\bf A}\}=\{\bd{x}\in\R^m|\bd{x}={\bf A}\theta, \quad\text{for } \theta\in\R^n,  \theta \ge {\bf 0}, \text{ and }{\bf 1}^\top\theta=1\}.
 \end{equation}	 
 Note that we will illustrate with ${\bf G}^\top={\bf P}^\top{\bf C}^\top$ in the geometric structure because gene expressions across sample issues, i.e.,  columns of ${\bf G}^\top$ (row of ${\bf G}$) are our major interested data features. By non-negativity of ${\bf C}$ and ${\bf P}$, we have
 	\begin{equation}
		{\bf G}_{(i)} \in {cone}({\bf P}^\top)\subseteq\R^{n}_+, \quad 1\le i \le N,
	\end{equation}
	and equivalently ${cone}({\bf G}^\top)\subseteq {cone}({\bf P}^\top)\subseteq\R^{n}_+$. Then problem (\ref{eqn:dconv}) can be interpreted as finding a nested cone problem: given two nested cones,  ${cone}({\bf G}^\top)$ and $\R^n_+$, find the nested cone ${cone}({\bf P}^\top)$ between them. 	
\begin{figure}
	\begin{center}
		\includegraphics[width=0.95\textwidth]{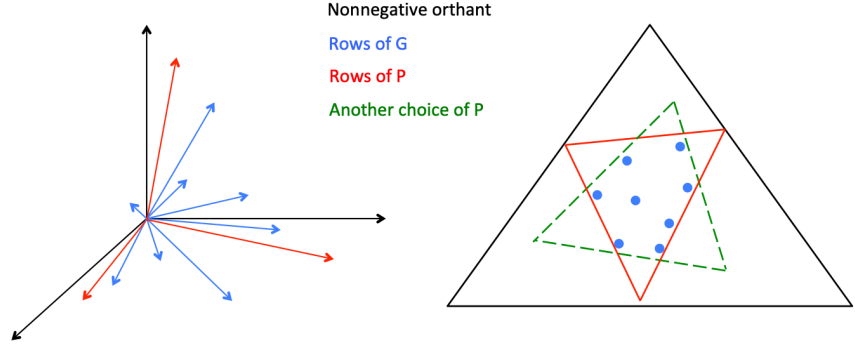} 
	\end{center}
	\caption{Cone (left) and convex hull (right) views of NMF as a nested cone problem.}%
\label{fig:geo}%
\end{figure}
This interpretation is displayed in the left panel of Figure \ref{fig:geo} with $N=8$ and $k=n=3$. It is easier to interpret this idea as shown in the right panel,  in terms of convex hull view, which is one dimension less than the cone view. This can be done easily with normalization of  columns of ${\bf G}^\top$ and ${\bf P}^\top$ to their unit $l_1$ norm.  From the right panel of Figure \ref{fig:geo}, we can also see that the solution of ${conv}({\bf P}^\top)$  is not unique: data $\bf G$ can be contained within two different cones (red solid and green dashed triangles) formed by different choices of matrix ${\bf P}$. It also motivates the idea that if rows of ${\bf G}$ ``spread out'' enough in the nonnegative orthant, ${cone}({\bf P}^\top)$ or  ${conv}({\bf P}^\top)$ may be unique. This condition turns into the following theories on the identifiability of the NMF.

\subsection{Strong and weak conditions on identifiability:}  There are two types of conditions for unique solution of problem (\ref{eqn:dconv}) when $\epsilon=0$. First, we need the following definitions for separable  and sufficient scattered matrices:

%\paragraph{Definition (Separable matrix)}
\begin{definition}
	The matrix ${\bf A}\in\R_+^{m\times n}$ is separable if ${cone}({\bf A})=\R^m_+$. 
\end{definition}

\begin{figure}
	\begin{center}
		\includegraphics[width=0.95\textwidth]{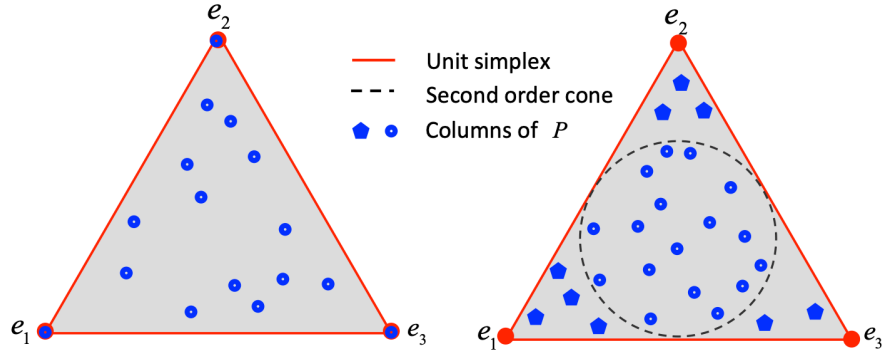} 
	\end{center}
	\caption{Convex hull views for strong (left) and weak (right) identifiability conditions with $k=3$.}%
\label{fig:ssc}%
\end{figure}

\begin{definition}\label{def:scattered}
	The matrix ${\bf A}\in\R_+^{m\times n}$ is sufficiently scattered if: (i) The second-order cone in $\R^m_+$ is contained in ${cone}({\bf A})$, i.e. $\mathcal{C}=\{{\bf x}\in\R^m_+|{\bf e}^\top {\bf x}\ge\sqrt{m-1}\|{\bf x}\|_2\}\subseteq {cone}({\bf A})$; {\em and} (ii) There does not exist any orthogonal matrix ${\bf Q}$ such that ${cone}({\bf A})\subseteq {cone}({\bf Q})$, except for permutation matrices.
\end{definition}
Based on these definitions, there are \cite{huang2013non,donoho2004does,laurberg2008theorems,gillis2011nonnegative}:
	\begin{theorem}[Strong identifiability condition]\label{thm:strong}
		Assuming $k = \text{rank}({\bf G})$, $\epsilon=0$, if problem (\ref{eqn:dconv}) admits a  solution, for which both ${\bf C}^\top$ and ${\bf P}$ are separable matrices, then the solution is unique. 
	\end{theorem}
	\begin{theorem}[Weak identifiability condition]\label{thm:weak}
		Assuming $k = \text{rank}({\bf G})$, $\epsilon=0$, if both ${\bf C}^\top$ and ${\bf P}$ are sufficiently scattered, then problem (\ref{eqn:dconv}) admits a unique solution. 
	\end{theorem}
	Theorem \ref{thm:strong} is quite rigorous: being separable matrices means ${\bf P}$ and ${\bf C}^\top$ must contains (scaled) extreme rays of the nonnegative orthant in the corresponding space, i.e., for every $r=1,2,...k$ there exists a column index $l_r$, such that ${\bf P}^{(l_r)}=\alpha_r{\bf e}_r\in\R^k_+$, where $\alpha_r$ is a scalar. This condition is illustrated in the left panel of Figure \ref{fig:ssc} as $k=3$. Note it is the convex hull view, so columns of ${\bf P}$  are represented as blue dots in the unit simplex (red triangle) in $\R^3_+$.  As in the figure, some columns of ${\bf P}$ are required to be exactly align with unit vectors ${\bf e}_r, r = 1, 2, 3$ (overlapping with red dots).  Similar situation is for matrix ${\bf C}^\top$. Such assumptions on both variables are too strong for practical applications, especially when noises present. 
	
	On the other hand, Theorem \ref{thm:weak}  is much more relaxed: The right panel of Figure \ref{fig:ssc} illustrates such condition: the dashed circle represents the intersection of the second-order cone in $\R_+^3$ and the unit simplex.  All blue objects, including dots and pentagons, are for columns of ${\bf P}$.  In this case, none of those columns are required to overlap with ${\bf e}_r$, but some of them (pentagons) need to fall out of the circle (second-order cone). 

\subsection{Relation to the gene expression data:} Figure \ref{fig:ssc} intuitively explained the strong and weak identifiability conditions of NMF  in order to achieve interpretable solutions. However, there are several issues when applying these theories to RNAseq data: 
(1) First of all and most importantly, how do these general identifiability conditions relate to the specific biological problem, i.e. bulk tissue RNAseq data?
(2) What measure should be used to define the ``sufficient scattering'' of matrix columns? Euclidean distance is used to illustrate the idea  in  Figure \ref{fig:ssc}, but it is not practical for high dimension data because column (row) normalization is needed for  convex hull description.  
(3) Theorems \ref{thm:strong}  and  \ref{thm:weak} can be used as  variable constraints in the optimization problem, but the actual questions are  how many, and which columns or (rows) the constraints should be enforced? How can we obtain this information from the only available data ${\bf G}$?  

The concept of  marker genes will help to address these issues.  By its name,  marker genes of a certain type of cell dominantly express in that cell type while rarely express in others. Each cell type may have multiple marker genes but one marker gene is only for one cell type. Mathematically, for each cell type $r = 1,2, 3,...k$, there exists an index set $\mathcal{S}_r$, such that for any $i\in\mathcal{S}_r$: expression level $c_{ir}$ is the dominant entry  (the only nonzero entry in ideal noiseless case) in the $i$-th row  of  matrix ${\bf C}$.  
This characteristics implies that ${\bf C}^\top$ is separable in ideal case and its columns are sufficiently scattered if noises present. On the other hand, no structures can be assumed for matrix ${\bf P}$. Although all conditions in Theorem \ref{thm:strong} or \ref{thm:weak} are not fully satisfied, reasonable solutions can be expected with constraints on variable ${\bf C}$.

Then the questions is how to identify marker genes (their index $\mathcal{S}_r$) from the data ${\bf G}$. Actually  in the noiseless case, for a given $r=1,2,..k$ and any $i\in \mathcal{S}_r$,  ${\bf C}_{(i)}=\alpha_i{\bf e}_r^\top$,  where ${\bf e}_r$ is the unit basis vector in $\R^r$. As consequences, the $i$-th row ${\bf G}_{(i)}=\alpha_i{\bf P}_{(r)}$ and hence all rows of ${\bf G}$ are linearly dependent if their indices are from the same set $\mathcal{S}_r$. In the more practical scenario where noise present, this linear dependence among vectors will become strong correlations.   Many literature has confirmed this phenomena that gene expressions across samples, or rows of ${\bf G}$, will display strong correlation, if they are from marker genes of the same cell type \cite{kuhn2012cell,mohammadi2016critical,zaitsev2019complete,avila2018computational,shen2010cell}.

\section{Mathematical Model:}\label{sec:model}

\begin{figure}[ptb]
\begin{center}%
\begin{tabular}[c]{cc}%
	\includegraphics[width=0.5\textwidth]{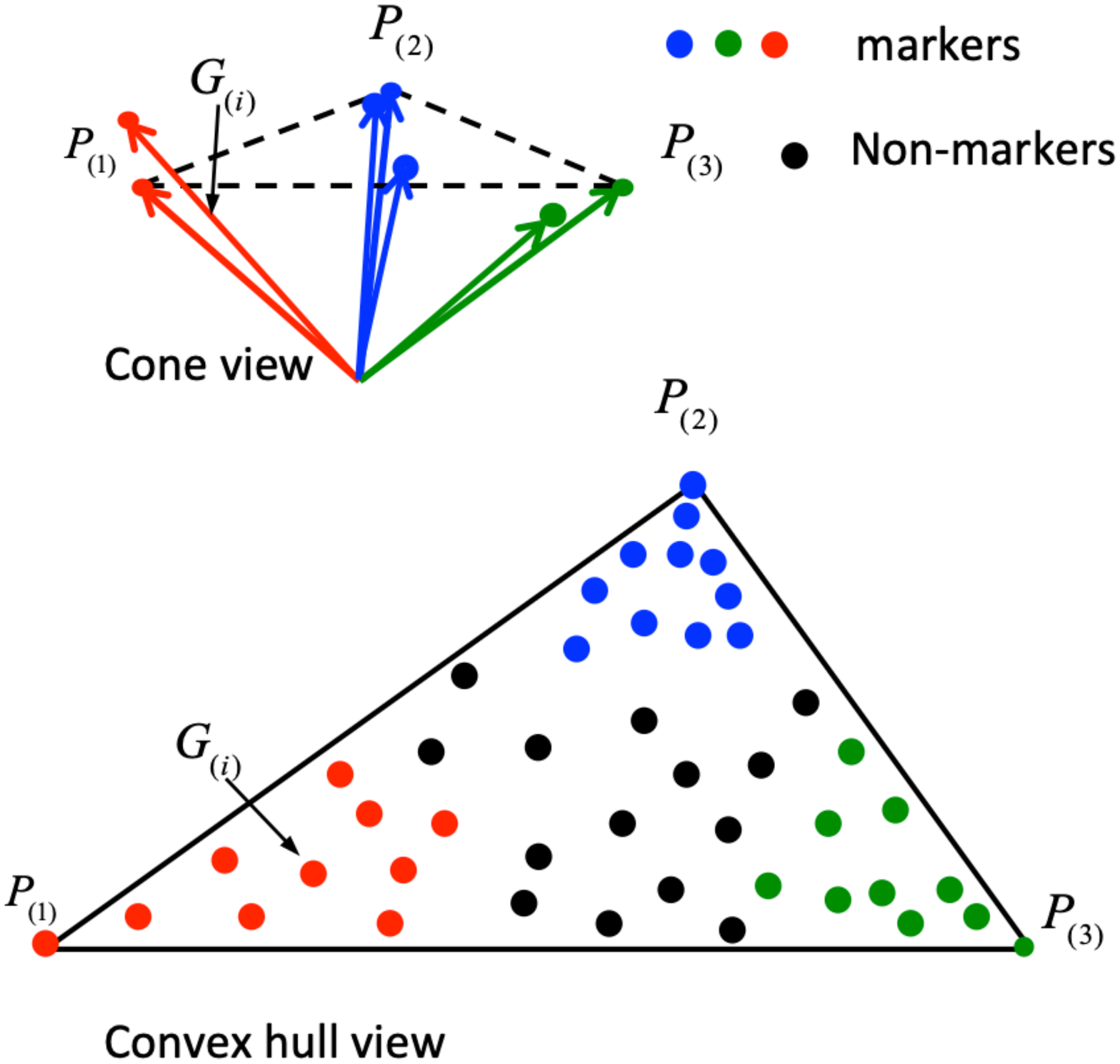} &
	\includegraphics[width=0.5\textwidth]{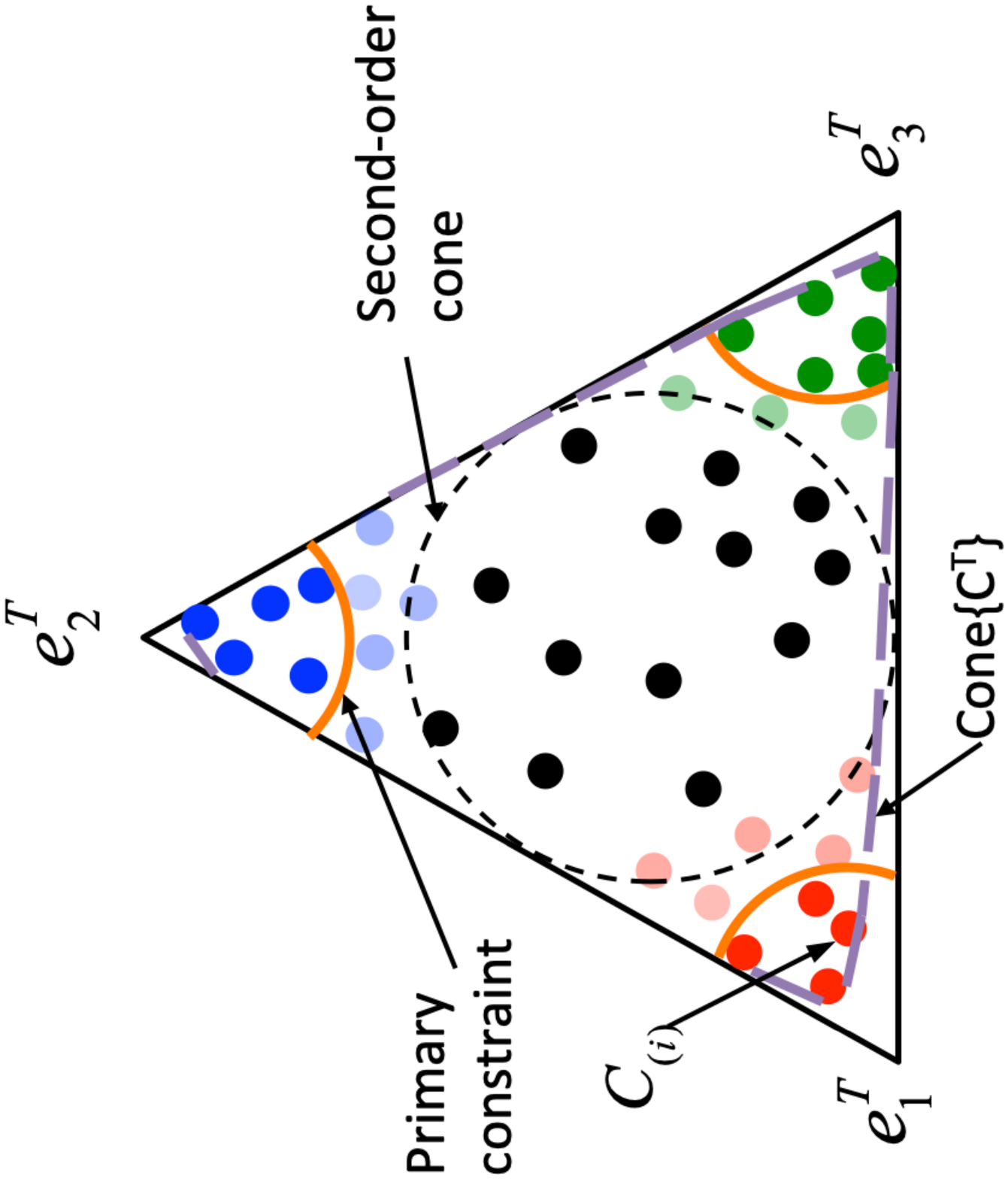}
	\end{tabular}
\end{center}
\caption{RNAseq data structure (left) and geometric constraints (right) }%
\label{fig:cluster-constraints}%
\end{figure}

Bulk tissue RNAseq data ${\bf G}$ has richer structures than just ${\bf G}_{(i)} \in {cone}({\bf P}^\top)\subseteq\R^{n}_+$:  for marker genes of the same cell type, e.g., the $r$-th cell type, their expressions across samples are highly correlated and correlated to ${\bf P}_{(r)}$. A cone view of this property is displayed in the left-top panel of Figure \ref{fig:cluster-constraints} for $k=3$. It is easier to investigate this feature further from the convex hull view, in which each row of ${\bf G}$ can be represented by a dot: As shown by the left-bottom panel of Figure \ref{fig:cluster-constraints},   rows of $\bf G$ for the marker genes of the  $r$-th cell type  tend to form a cluster around ${\bf P}_{(r)}$ due to the strong correlation.  This property motivates us that it is possible to identify marker genes from data $\bf G$ by clustering its rows and to quantitatively explore the  geometric structures of its row space. Further, note that ${\bf C}_{(i)}$ is actually the coefficient vector of ${\bf G}_{(i)}$ under the basis vectors of rows of ${\bf P}$, we can transfer such geometric structure of ${\bf G}$ to ${\bf C}$, hence to enforce the weak identifiability condition on variable $\bf C$. The two steps are termed as finding and preserving the structure, respectively, which will be detailed as the following:

\subsection{Finding geometric structures by spectral clustering analysis}

In this step we will classify ${\bf G}_{(i)}, 1\le i \le N$ into $k$ groups to identify possible marker genes for the $k$ types of cells.  Among many existing clustering techniques, we propose to use spectral clustering \cite{von2007tutorial}, which is one type of manifold learning algorithms that can explore intrinsic geometric/topological structure of high dimensional data. Thus, it has many fundamental advantages and very often outperforms traditional clustering algorithms such as $k$-means or single linkage. 
	 To perform spectral clustering, we need the similarity graph $G=(V,E)$, with vertex set $V=\{{\bf G}_{(i)}\}_{i=1}^N\subset\R^{n}$. The non-negative weights $\omega_{ij}$ of edges $E=\{e_{ij}\}$ are calculated by a   function $\R^n\times\R^n\to\R_+$, quantifying the correlation between two vertices.  
	 %Vertices of $G$, or row vectors of ${\bf G}$, will be clustered into groups by classic clustering algorithms applying to the relevant eigenvectors of the graph Laplacians.
	We propose to evaluate $\omega_{ij}$ as
 	\begin{equation}\label{eqn:omga}
		\omega_{ij}=\exp\left\{-\frac{d_{eisen}\left({\bf G}_{(i)}, {\bf G}_{(j)}\right)^2}{\sigma}\right\}, 1\le i \le N, 1\le j \le N,
	\end{equation}
	where
	\begin{equation}\label{eqn:eisen}
		d_{eisen}\left({\bf G}_{(i)}, {\bf G}_{(j)}\right)=1-\frac{<{\bf G}_{(i)}, {\bf G}_{(j)}>}{|{\bf G}_{(i)}||{\bf G}_{(j)}|}
	\end{equation}
	is the Eisen cosine correlation distance and $\sigma>0$ is a parameter.
	The matrix ${\bf W}=(\omega_{ij})\in\R^{N\times N}$ is called the adjacency matrix. Meanwhile, define its 
degree matrix  ${\bf D}={\rm diag}(d_1,d_2,...d_N)$, where $d_i=\sum_{j=1}^N\omega_{ij}$ is the degree of the vertex ${\bf G}_{(i)}$.
	 With these matrices,  different types of graph Laplacians (gL) of  $G$ can be defined, such as the unnormalized gL $ L={\bf D-W}$, symmetric normalized  gL $L_{\rm sym}={\bf I}-{\bf D}^{-\frac{1}{2}}{\bf W}{\bf D}^{-\frac{1}{2}}$,  or random walk gL $L_{\rm rw}={\bf I}-{\bf D}^{-1}{\bf W}$. By examining the first a few eigenvectors of gL,  rows of data $\bf G$ will be clustered into $k$ groups and the sets $\{\mathcal{G}_r\}_{r=1}^k$ record row indices of ${\bf G}$ in the corresponding clusters. Choice of different gLs depends on specific data applications \cite{qin2019fast}. For our problem, we use the normalized gL $L_{\rm sym}$ and  perform the spectral clustering package in  Matlab. 
	 
	 %Further,  we  need to find a subset for each $\mathcal{G}_r$, i.e, $\mathcal{S}_r\subset\mathcal{G}_r$, to represent the indices of marker genes for each cell type.

%Add introduction and brief process of spectral clustering here.

\subsection{Geometric structure guided model} 

With the clustering information, we are able to establish geometric structure guided model by applying two constraints  on the row space of variable $\bf C$: the solvability constraint and manifold regularization. 
This work is motivated by manifold regularization which uses graph Laplacian as regularizer, but it  provides new characteristics. The major novelty is to incorporate the identifiability conditions  in the regularization of the NMF (primary constraint). Another new feature is to encode the geometric information of the data space   (secondary constraint)  based on Eisen cosine correlation distance, instead of Euclidean  distance in traditional graph regularized NMF.

\subsection{Solvability constraint: } 
According to Theorem \ref{thm:weak}, rows of ${\bf C}$ need to be scattered sufficiently, such that the second-order cone in $\R_+^r$ is contained in ${cone}({\bf C}^\top)$.  Further, Definition \ref{def:scattered} implies that this requirement is only needed for {\em some} of rows in $\bf C$, i.e., those rows corresponding to marker genes. 
In the previous step, all rows of ${\bf G}$ are clustered into $k$ groups and their row numbers are recorded in the set $\{\mathcal{G}_r\}_{r=1}^k$. In the current step,  correlations are ranked within each group, and a subset, i.e., $\mathcal{S}_r\subset\mathcal{G}_r$ is determined accordingly to represent the indices of marker genes for each cell type. Note that row indices of both ${\bf G}$ and ${\bf C}$ represent gene IDs, we need rows of $\bf C$ with index $\mathcal{S}_r$ to scatter enough to accommodate the second-order cone. This can be done
by requiring them to have strong correlations with $\{{\bf e}_r^\top\}$, i.e.,  defining the penalty function: 
	\begin{equation}\label{eqn:con1}
		%\mathcal{F}_1({\bf C})=\sum_{r=1}^k\sum_{i\in\mathcal{S}_r}\left(|{\bf C}_{(i)}||{\bf e}^\top_{r}|-<{\bf C}_{(i)}, {\bf e}^\top_{r}>\right).
		\mathcal{F}_1({\bf C})=\frac{\lambda_1}{2}\sum_{r=1}^k\sum_{i\in\mathcal{S}_r}d_{eisen}\left({\bf C}_{(i)}, {\bf e}^\top_{r}\right)^2,
	\end{equation}
	where $\lambda_1$ is a parameter. 
	 This idea is illustrated from the convex hull view   in the right panel of Figure \ref{fig:cluster-constraints} as $k=3$. Red, blue, green dots represent rows of $\bf C$ that indexed by $\mathcal{G}_r$. The darker colored dots, representing selected marker genes in $\mathcal{S}_r$, are ``required'' to stay in the circular sectors (orange), such that the ${conv}\{\bf C^\top\}$ (dashed purple) is large enough to contain the second-order cone (dashed circle). With the Eisen cosine correlation distance, we do not need to work in ${conv}\{{\bf C}^\top\}$ (which requires normalizing coefficients), but directly in ${cone}\{{\bf C}^\top\}$. 	
	 
		%Note that the formula in the parenthesis   has the same effect as the formula $1-<{\bf C}_{(i)}, {\bf e}^\top_{r}>/|{\bf C}_{(i)}||{\bf e}^\top_{r}|$ to force the angles small, but Eq. (\ref{eqn:con1}) will give much simpler formula for the gradient, as shown in the next section. Additionally, when it is rigorously zero, Eq. (\ref{eqn:con1}) represents the strong separability condition.

	  \subsection{Manifold constraint: } 
	  According to   the local invariance assumption in manifold regularization \cite{belkin2002laplacian,cai2009probabilistic,he2004locality}, if two data points  are close in the intrinsic geometry of the data distribution, then the representations of this two points in a new basis  should be also close to each other under the same metric.
	  Note that ${\bf C}_{(i)}$ is the representation of the data point ${\bf G}_{(i)}$ under the basis ${\bf P}^\top$, then by	 such manifold assumption, we 	require matrix $\bf C$ to inherit the similar geometric structure of matrix ${\bf G}$, i.e., rows of $\bf C$ belong to the same cluster have strong mutual correlations. To achieve this goal, we define another penalty function:
	  \begin{equation}\label{eqn:con2}
		\mathcal{F}_2({\bf C})=\frac{\lambda_2}{2}\sum_{j=1}^N\sum_{i=1}^N\omega_{ij}d_{eisen}\left({\bf C}_{(i)}, {\bf C}_{(j)}\right)^2.
	\end{equation}
	Recall that  entry $\omega_{ij}>0$ in the adjacency matrix ${\bf W}$ in (\ref{eqn:omga})  measure the correlations (larger value  represents stronger correlation)  between genes $i$ and $j$ in data $\bf G$.
	
	\vspace{10pt}
	\noindent{\em Remark 1:} Equation (\ref{eqn:con2})  is a generalization of  the traditional graph Laplacian regularization \cite{cai2010graph, qin2019fast}. Actually,  if the Eisen cosine correlation distance in  (\ref{eqn:con2}) is replaced by the Euclidean distance ($l_2$ norm), then
	 \begin{eqnarray}\nonumber
		\mathcal{F}_2({\bf C})&=&\frac{\lambda_2}{2}\sum_{j=1}^N\sum_{i=1}^N\omega_{ij}||{\bf C}_{(i)}- {\bf C}_{(j)}||^2
		=\sum_{i=1}^Nd_{ii}{\bf C}^\top_{(i)}{\bf C}_{(i)}-\sum_{i=1}^N\omega_{ij}\sum_{j=1}^N{\bf C}^\top_{(i)}{\bf C}_{(j)}\\\nonumber
		&=&\lambda_2{\rm Tr}({\bf C}^\top L{\bf C}),
	\end{eqnarray}
	with $L$ being the graph Laplacian operator defined earlier.
	 
	 \vspace{10pt}
	 \noindent{\em Remark 2:} The current work is more than a generalization of traditional graph Laplacian regularized NMF by using correlation distance metric. Indeed, the manifold assumption only requires ${\bf C}_{(i)}$ and ${\bf C}_{(j)}$ to be close, if ${\bf G}_{(i)}$ and ${\bf G}_{(j)}$ are close in the distance metric, i.e., genes $i$ and $j$ are classified to belong the {\em same} cell type. However, it does not require rows of ${\bf C}$ to be far from each other if they are for {\em different} types of cells. This issue is addressed by Eq. (\ref{eqn:con1}) and the novelty of the proposed work is inclusion of the solvability constraints.
	 
	  	 \vspace{10pt}
		 \noindent {\em Remark 3:} In the extreme case $\lambda_1\to\infty$, we have ${\bf C}_{(i)}=\alpha_i{\bf e}^\top_{r}$ if $i\in\mathcal{S}_r$. This case corresponds to the strong identifiability condition. However, in the realistic circumstances where noises present, this extreme requirement does not provide the optimal results, as shown in  numerical simulations. 
			
	\vspace{10pt}
	\noindent{\em Remark 4:}
	Relations between constraints (\ref{eqn:con1}) and (\ref{eqn:con2}) can be explained by the right panel of Figure \ref{fig:cluster-constraints} as $k=3$.  All genes  are classified into $k$ groups (red, blue, and green), indexed by $\mathcal{G}_r$,  assuming there are $k$ types of cells in those tissue samples.  In addition to the constraints that darker colored dots (indexed by $\mathcal{S}_r$) close to ${\bf e}_r^\top$, the dots in the same color need to stay close, corresponding to the geometric structure of data $\bf G$ shown in the left panel.

\subsection{Full model:}Combining the solvability condition (\ref{eqn:con1}) and manifold constraint (\ref{eqn:con2}),
	we propose the following geometric structure guided nonnegative matrix factorization (GS-NMF) model. For illustration convenience, we define the set $T:=\{Z\in\R_+^{k\times n}, {\bf 1}^\top Z={\bf 1}^\top\}$
	and the indicator function $\mathds{1}_{T}$  as $\mathds{1}_T(Z) = 0$ if      $Z\in T$ while  $\mathds{1}_T(Z) = \infty$ otherwise.
%	 \begin{equation}
%	 	 \mathds{1}_G(Z) = 
%		\begin{cases}
%    				0       & \quad \text{if } Z\in G\\
%    				\infty  & \quad \text{otherwise} 
% 		 \end{cases}
%	 \end{equation}
	With these notations, solving for $\bf C$ and $\bf P$ becomes the optimization problem: 
	\begin{equation}\label{eqn:model}
		\min_{{\bf C}\ge 0, {\bf P}\ge 0}{\frac{1}{2}\|{\bf G}-{\bf C}{\bf P}\|_F^2+\mathcal{F}({\bf C})+ \mathds{1}_T({\bf P}). }
	\end{equation}
	%where $\mathcal{F}({\bf C})=\mu\mathcal{F}_1({\bf C}) + (1-\mu)\mathcal{F}_2({\bf C})$ for some $0\le \mu\le 1$.
	In the first term, Frobenius norm to measure the error between deconvoluted solution and the given data. The total regularization function
	$\mathcal{F}({\bf C})= \mathcal{F}_1({\bf C}) +  \mathcal{F}_2({\bf C})$ as each component is defined in (\ref{eqn:con1})-(\ref{eqn:con2}). The third term $\mathds{1}_T({\bf P})$  simply means  sum-to-one conditions on {\em columns} of $\bf P$, or column stochasticity, since the sum of cellular proportions in each tissue sample is supposed to be one.

	\subsection{Gradient: } For computational algorithms, it is useful to derive the gradients of Eq. (\ref{eqn:con1}) and (\ref{eqn:con2}).
		A detailed element-wise computation of $\displaystyle{\frac{\partial \mathcal{F}_1}{\partial c_{ij}}}$ yields	
	\begin{equation}\label{eqn:grad1}
		\frac{\partial\mathcal{F}_1}{\partial c_{ij}}=\lambda_1\chi_{g(i)}\left(1-\frac{c_{ig(i)}}{|{\bf C}_{(i)}|}\right)\frac{1}{|{\bf C}_{(i)}|}\left[\frac{c_{ig(i)}}{|{\bf C}_{(i)}|^2}c_{ij}-\delta_{j,g(i)}\right],
			\end{equation}
	where $g:\mathcal{S}_r\to r$ is a surjective function, i.e., $g(i)=r$ if there is an $r$ such that $i\in\mathcal{S}_{r}$ or $g(i)=0$ otherwise. This function is known from the spectral clustering and $g(i)=r$ means gene $i$ is the marker gene of the $r$-th cell type. To take into account all possible non-marker genes, we use characteristic function $\chi_{g(i)}=1$ if $g(i)\neq 0$ and $\chi_{g(i)}=0$ otherwise. It is convenient to write Eq. (\ref{eqn:grad1}) in matrix form for future computation. To do so, we define
	 ${\bf W}_1={\rm diag}\displaystyle{\left\{\left(1-\frac{c_{ig(i)}}{|{\bf C}_{(i)}|}\right)\frac{\chi_{g(i)}}{|{\bf C}_{(i)}|}\right\}}$;  ${\bf W}_2={\rm diag}\displaystyle{\left\{\frac{c_{ig(i)}}{|{\bf C}_{(i)}|^2}\right\}}$, and ${\bf C}_g=\{\delta_{g(i), j}\}$. 
	 Here we use the notation that ${\rm diag}\{x_i\}$ means a diagonal matrix with all entries $x_i$ on its diagonal and ${\rm diag}({\bf x})$ means a vector $\bf x$ forms a diagonal matrix with the corresponding size. Then the matrix form of Eq. (\ref{eqn:grad1}) is
	\begin{equation} \label{eqn:mgrad1}
	 	\nabla\mathcal{F}_1 =\lambda_1{\bf W}_1({\bf W}_2{\bf C}-{\bf C}_g).
	\end{equation}
	
	For the second constraint (\ref{eqn:con2}), the gradient is
	\begin{eqnarray}\nonumber
		\frac{\partial\mathcal{F}_2}{\partial c_{ij}}&=&\lambda_2\left[\sum_{l=1}^N\omega_{il}\left(1-\frac{<{\bf C}_{(i)}, {\bf C}_{(l)}>}{|{\bf C}_{(i)}||{\bf C}_{(l)}|}\right)\frac{<{\bf C}_{(i)}, {\bf C}_{(l)}>}{|{\bf C}_{(i)}||{\bf C}_{(l)}|}\right]\frac{c_{ij}}{|{\bf C}_{(i)}|^2}\\\nonumber
		&-&2\lambda_2\sum_{l=1}^N\omega_{il}\left(1-\frac{<{\bf C}_{(i)}, {\bf C}_{(l)}>}{|{\bf C}_{(i)}||{\bf C}_{(l)}|}\right)\frac{1}{|{\bf C}_{(i)}||{\bf C}_{(l)}|}c_{lj}.
	\end{eqnarray}

	To have matrix formulation, we define the   matrix ${\bf coC} = \left\{\frac{<{\bf C}_{(i)}, {\bf C}_{(j)}>}{|{\bf C}_{(i)}||{\bf C}_{(j)}|}\right\}$,  the column vector $|{\bf C}|$  with entries being $l_2$ norms of rows of ${\bf C}$, and $|{\bf C}|^{-1}$ as its element-wise reciprocal.  Further, define ${\bf W}_3={\bf W}\odot({\bf 1}_{N\times N}-{\bf coC})\odot {\bf coC}$, and  ${\bf W}_4={\bf W}\odot({\bf 1}_{N\times N}-{\bf coC})\odot (|{\bf C}|^{-1}(|{\bf C}|^{-1})^{\top}$), with ``$\odot$'' being the Hadamard  product of matrices. Consequently, the matrix form of the gradient of Eq. (\ref{eqn:con2}) is
	\begin{equation}\label{eqn:mgrad2}
		%\nabla\mathcal{F}_2 =4\lambda_2\left[{\rm sum}({\bf W}_1, 2){\rm diag}(|{\bf A}|^{-2})-{\bf W}_2\right]{\bf A},
		\nabla\mathcal{F}_2 = \lambda_2{\rm diag}({\bf W}_3{\bf 1}){\rm diag}(|{\bf C}|^{-2}){\bf C} - 2\lambda_2{\bf W}_4{\bf C}.
	\end{equation}

	%Assembling the geometric constraint to ${\bf C}$ and column stochastic constraint to ${\bf P}$ into the NMF, 

	 \section{Computational  algorithms} \label{sec:algorithm}
	 
	 It is well-known that the objective function in (\ref{eqn:model}) is non-convex in both variables together. There are several types of numerical methods to obtain a local minimum, including Multiplicative Update Algorithm (MUA), Alternating nonnegativity constrained least squares (ANLS), and  the alternating direction method of multipliers (ADMM), etc. The first two types of methods are quite straightforward for the basic NMF problems, or with linear constraints. While for the nonlinear constraints (\ref{eqn:mgrad1})-(\ref{eqn:mgrad2}), it is convenient to adopt ADMM framework to develop numerical schemes.

	 To do so, we first rewrite model (\ref{eqn:model}) as
	 \begin{equation}\label{eqn:model2}
		\min_{{\bf C}, {\bf P}}{\frac{1}{2}\|{\bf G}-{\bf C}{\bf P}\|_F^2+\mathcal{F}({\bf C})+\mathds{1}_S({\bf C}) + \mathds{ 1}_T({\bf P}) },
	\end{equation}
	where $S:=\R_+^{N\times k}$ is the set of all nonnegative matrices of the size $N\times k$. Then we introduce two auxiliary variables ${\bf A}$ and $\bf Q$, and rewrite (\ref{eqn:model2}) into an equivalent form
	\begin{equation}
		\begin{aligned}
		& \min_{{\bf C}, {\bf P}, {\bf A}, {\bf Q}}{\frac{1}{2}\|{\bf G}-{\bf C}{\bf P}\|_F^2}+\mathcal{F}({\bf A})+\mathds{1}_S({\bf C}) + \mathds{1}_T({\bf Q}) \\
		& \text{s.t.: } {\bf C}-{\bf A} = {\bf 0}, {\bf P}-{\bf Q} = {\bf 0},
		\end{aligned}
	\end{equation}
	and the corresponding augmented Lagrange function is \cite{warren2015hyperspectral} 
	\begin{equation}\label{eqn:admm}
		\mathcal{L} = \frac{1}{2}\|{\bf G}-{\bf C}{\bf P}\|_F^2+\mathcal{F}({\bf A})+\mathds{1}_S({\bf C}) + \mathds {1}_T({\bf Q})
			   	  + \frac{\rho}{2}\|{\bf C}-{\bf A}+\tilde{\bf A}\|_F^2+\frac{\gamma}{2}\|{\bf P}-{\bf Q} + \tilde{\bf Q}\|_F^2
	\end{equation}
	where $\tilde{A}$ and $\tilde{Q}$ are dual variables of ${\bf A}$ and $\bf Q$, respectively, while $\rho>0$ and $\gamma>0$ are penalty parameters. As results, the ADMM (\ref{eqn:admm}) can be written as an iteration (from $i$- to $i+1$-th step) in its scaled form \cite{boyd2011distributed}:
	\begin{equation}\label{eqn:scaledADMM}
		\begin{aligned}
		& {\bf C}^{i+1}:= \argmin_{{\bf C}}{\frac{1}{2}\|{\bf G}-{\bf C}{\bf P}^i\|_F^2+\frac{\rho}{2}\|{\bf C}-{\bf A}^i+\tilde{\bf A}^i\|_F^2+\mathds{1}_S({\bf C})}\\
		& {\bf P}^{i+1}:= \argmin_{{\bf P}}{\frac{1}{2}\|{\bf G}-{\bf C}^i{\bf P}\|_F^2+\frac{\gamma}{2}\|{\bf P}-{\bf Q}^i+\tilde{\bf Q}^i\|_F^2+ \mathds{1}_T({\bf P}) }\\
		& {\bf A}^{i+1}:= \argmin_{{\bf A}}{\mathcal{F}({\bf A})+\frac{\rho}{2}\|{\bf C}^i-{\bf A}+\tilde{\bf A}^i\|_F^2}\\
		& {\bf Q}^{i+1}:= \argmin_{{\bf Q}}{\frac{\gamma}{2}\|{\bf P}^i-{\bf Q}+\tilde{\bf Q}^i\|_F^2}\\
		& \tilde{\bf A}^{i+1}:=\tilde{\bf A}^i+{\bf C}^i-{\bf A}^i\\
		& \tilde{\bf Q}^{i+1}:=\tilde{\bf Q}^i+{\bf P}^i-{\bf Q}^i
		\end{aligned}.
	\end{equation}
	Each variable in system (\ref{eqn:scaledADMM})  can be solved individually. Specifically,
	for the ${\bf C}$-subproblem, the Karush-Kuhn-Tucker (KKT) condition \cite{nocedal2006numerical} yields a closed form for ${\bf C}$, i.e.,
	\begin{equation}\label{eqn:c}
		{\bf C}^{i+1}=[{\bf G}{\bf P}^{i^\top}+\rho({\bf A}^i-\tilde{\bf A}^i)]({\bf P}^{i}{\bf P}^{i^\top}+\rho{\bf I})^{-1},
	\end{equation}
	where ${\bf P}^{i}{\bf P}^{i^\top}+\rho{\bf I}$ is a small $k\times k$ matrix that can be inverted easily. The nonnegativity of ${\bf C}$ is obtained by row-wise active set method.
	%use active set method, look at the code in matlab to see how to set active an passive sets.
%	Can be solved with nonnegative LS algorithms row-wisely.
%	\begin{equation}
%		\begin{aligned}
%		& \min_{\bar{\bf C}}{\frac{1}{2}\sum_{j=1}^m\|{\bf M}\bar{\bf C}_{(j)}-{\bf b}_j\|_2^2} \\
%		& \text{s.t.: } \bar{\bf C}_{(j)}\ge 0, j=1,2,...m,
%		\end{aligned}
%	\end{equation}
%	where ${\bf M}={\bf P}^{i^T}+\rho{\bf I}$ and ${\bf b}_j=\bar{\bf X}_{(j)}+\rho\left({\bf A}^i_{(j)}-\tilde{\bf A}^i_{(j)}\right)$
%
	For the ${\bf P}$-subproblem,  KKT condition gives
	\begin{equation}\label{eqn:p}
		{\bf P}^{i+1}=\Pi\{({\bf C}^{i^\top}{\bf C}^i+\gamma{\bf I})^{-1}[{\bf C}^{i^\top}{\bf G}+\gamma({\bf Q}^i-\tilde{\bf Q}^i)]\},
	\end{equation}
	and this is a small-scale problem, in which a $k\times k$ matrix is to be inverted with column-wise  probability simplex projection $\Pi$ \cite{wang2013projection}.
	Solution of the ${\bf Q}$-subproblem is simply
	\begin{equation}\label{eqn:q}
		{\bf Q}^{i+1}=\max\{{\bf P}^i+\tilde{\bf Q}^i, {\bf 0}\}.
	\end{equation}

	The ${\bf A}$-subproblem involves the solvability condition (\ref{eqn:con1}) and manifold constraints (\ref{eqn:con2}), both of which are non-linear problem so no closed form can be used. To solve this subproblem we have to use the gradient descent method and make this step an inner iteration.
%	\begin{equation}
%		f({\bf A})=\frac{\lambda}{2}\sum_{i,j=1}^m\tilde{W}_{ij}\left(|{\bf A}_{(i)}||{\bf A}_{(j)}|-<{\bf A}_{(i)}, {\bf A}_{(j)}>\right)+\frac{\rho}{2}\|{\bf C}-{\bf A}+\tilde{\bf A}\|_F^2
%	\end{equation}
Denote the total objective function of the ${\bf A}$-subproblem as
	\begin{equation}
		f({\bf A})=\mathcal{F}({\bf A})+\frac{\rho}{2}\|{\bf C}-{\bf A}+\tilde{\bf A}\|_F^2,
	\end{equation}
	then according to (\ref{eqn:mgrad1}) and (\ref{eqn:mgrad2}),  its gradient is
	\begin{equation}\label{eqn:gA}
		\nabla f = \lambda_1{\bf W}_1({\bf W}_2{\bf A}-{\bf C}_g) +\lambda_2{\rm diag}({\bf W}_3{\bf 1}){\rm diag}(|{\bf C}|^{-2}){\bf A} - 2\lambda_2{\bf W}_4{\bf A} + \rho({\bf A-C}-\tilde{\bf A}).
	\end{equation}
	Note that Eq. (\ref{eqn:gA}) is nonlinear in terms of ${\bf A}$. For computational efficiency, we will use the result of ${\bf C}$ in the current outer iteration step to compute ${\bf W}_1$, ${\bf W}_2$, ${\bf W}_3$ and ${\bf W}_4$ so their values will not update in the inner loop.
	Algorithm \ref{alg:NMF} summarizes the entire processes of the GS-NMF. Note that  necessary raw data processing and spectral clustering steps are not included. The stoping criteria are to set $\|{\bf C}^{i+1}-{\bf C}^{i}\|/\|{\bf C}^i\|_F$ and $\|{\bf P}^{i+1}-{\bf P}^{i}\|/\|{\bf P}^i\|_F$ smaller than some tolerance.
	\begin{algorithm}%\label{alg:RPGD}
		\caption{Geometric structure constrained NMF (GS-NMF)}\label{alg:NMF}
		\begin{algorithmic}[1]
			\Require Data ${\bf G}$, initial guesses ${\bf C}_0$, ${\bf P}_0$,  structure identifier ${\bf C}_g$, graph adjacency matrix ${\bf W}$, tolerance $\epsilon$, parameters $\lambda_1$, $\lambda_2$, $\rho$ and $\gamma$.
			\Ensure   Matrices $\bf C$ and $\bf P$.
			\For{$i=0,1,...$ until criteria is satisfied} \quad\quad\quad\quad\quad\quad\quad\% outer iteration
					\State Solve ${\bf C}$-subproblem in Eq. (\ref{eqn:scaledADMM}) by (\ref{eqn:c});
					\State Solve ${\bf P}$-subproblem in Eq. (\ref{eqn:scaledADMM}) by (\ref{eqn:p});
					\For{$m=0,1,...$ until criteria is satisfied} \quad\quad\quad\quad\quad\quad\quad\% inner iteration
						\State Solve ${\bf A}$-subproblem in Eq. (\ref{eqn:scaledADMM}) through gradient descent method with gradient (\ref{eqn:gA});
					\EndFor
					\State Solve ${\bf Q}$-subproblem in Eq. (\ref{eqn:scaledADMM}) by (\ref{eqn:q});
					\State Set $\tilde{\bf A}^{i+1}:=\tilde{\bf A}^i+{\bf C}^i-{\bf A}^i$;
					\State Set $\tilde{\bf Q}^{i+1}:=\tilde{\bf Q}^i+{\bf P}^i-{\bf Q}^i$.
			\EndFor
		\end{algorithmic}
	\end{algorithm}

	%The reason we Do NOT assume sparseness of C is for potential noise quantification. Find C and then compare with its ``sparse form'', to derive error term.

%\begin{itemize}
%	\item Challenges (ill-posedness of the problem) and strategies (explore/preserve topological structure by finding marker genes): use marker gene to find $P$ and part of $C$, and then use fixed $P$ to find the rest part of $C$ 
%	\item Use spectral clustering to identify marker genes;
%	\item Model 1: Blind separation with Graph Laplacian regularization model.
%	\item Model 2:  manifold based separation for solving ${\bf C}$;
%	\item Algorithm1: conventional ADMM, Algorithm 2: manifold ADMM
%\end{itemize}

\section{Numerical Results}\label{sec:num}

In this section, we test   the proposed GS-NMF algorithms on two types of data.

\subsection{Simulations on synthetic data}

In order to have flexibility of matrix dimensions, noise levels, and known ground truth, we first test the algorithms on synthetic data, which are generated as the following strategies: matrix ${\bf C}\in\R_+^{N\times k}$ has a structure, so first we split $N = N_1+N_2+...N_k + N_{k+1}$. In order to mimic the marker gene expression in the corresponding cells, we generate $N_l, l = 1,2,...k$ rows of ${\bf C}$ such that they have strong correlations to ${\bf e}_l^\top$. The rest of $N_{k+1}$ rows are non-marker genes, so they are just generated randomly. Then all rows of ${\bf C}$ are  assembled and a random row-permutation is performed. Matrix ${\bf P}$ is just a  random ${k\times n}$ matrix with non-negative entries and columns normalized by their $l_1$ norms. Data ${\bf G}$ is computed simply by ${\bf G} = {\bf CP}+\epsilon$, where $\epsilon$ is the noise matrix following normal distribution.

\begin{figure}[h!]
\begin{center}%
\begin{tabular}[c]{cc}%
	\includegraphics[width=0.45\textwidth]{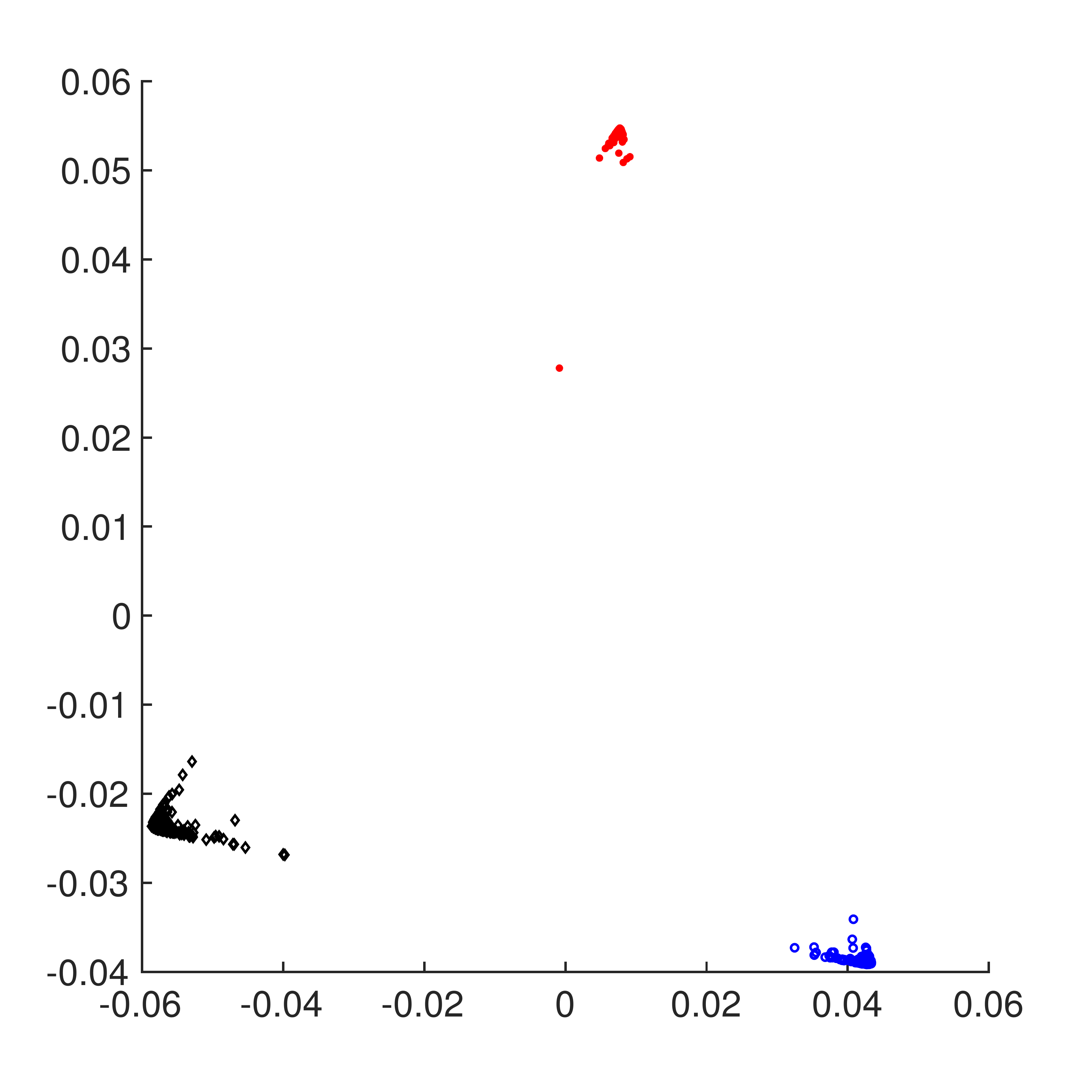} &
	\includegraphics[width=0.45\textwidth]{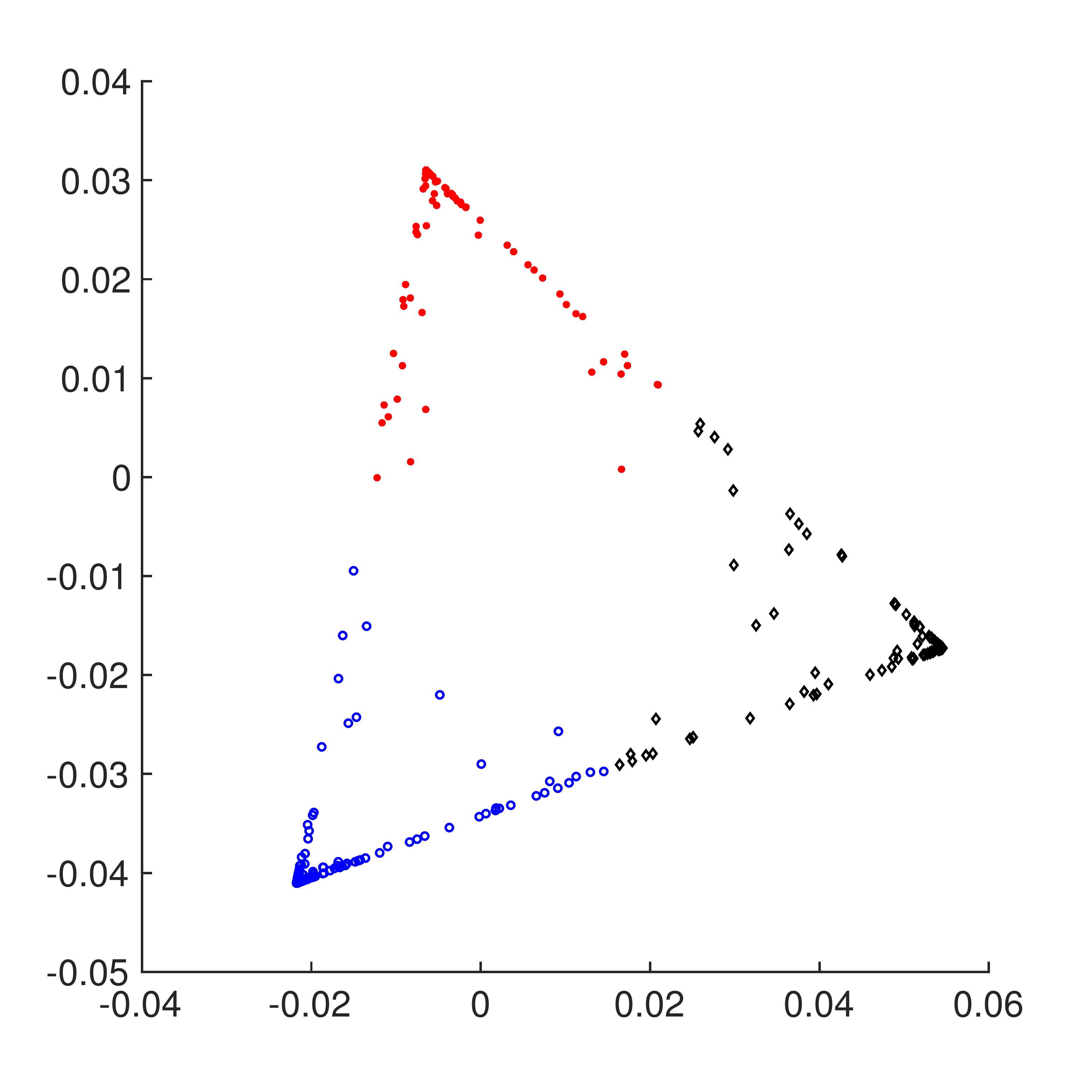}\\
	(a) & (b)\\
	\includegraphics[width=0.45\textwidth]{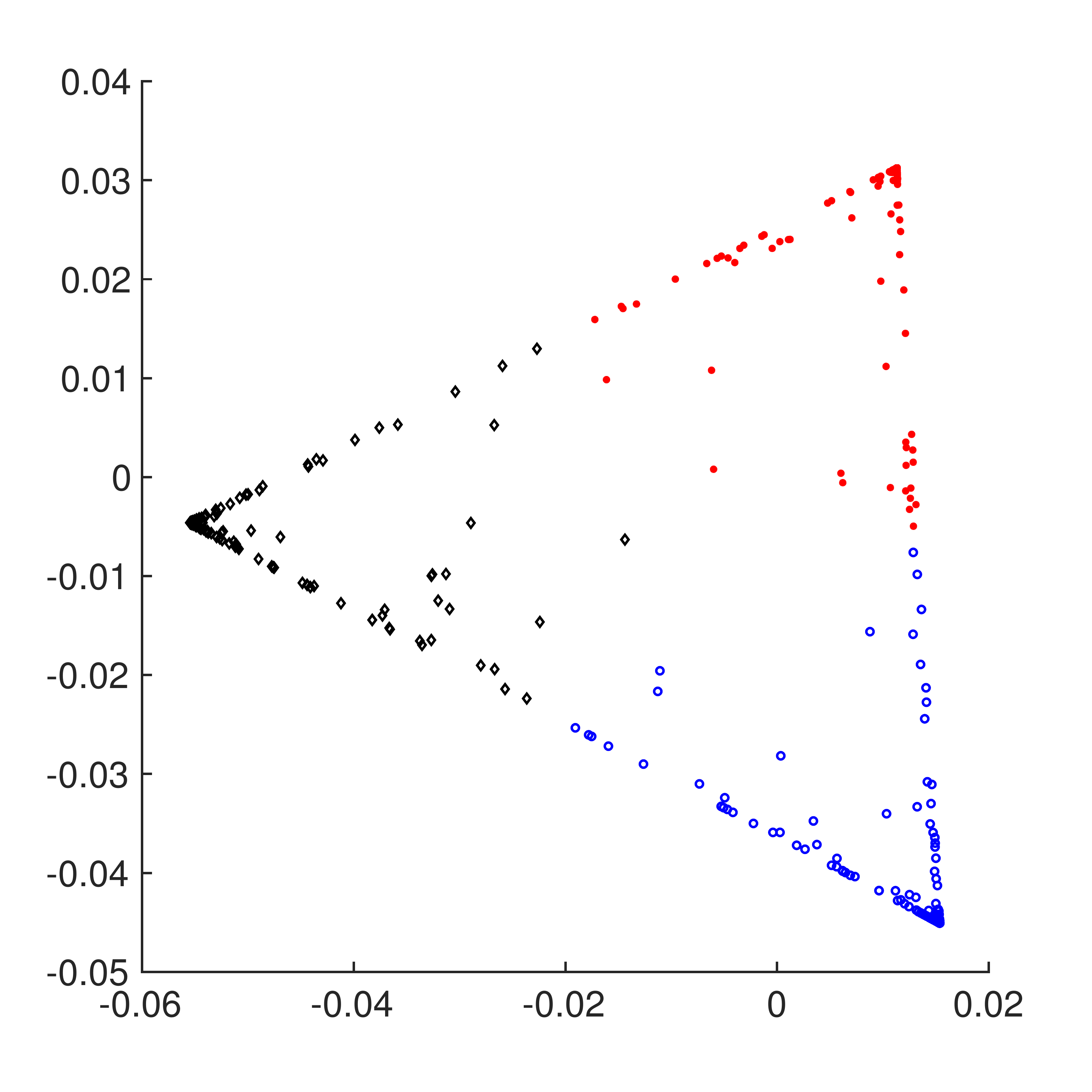} &
	\includegraphics[width=0.45\textwidth]{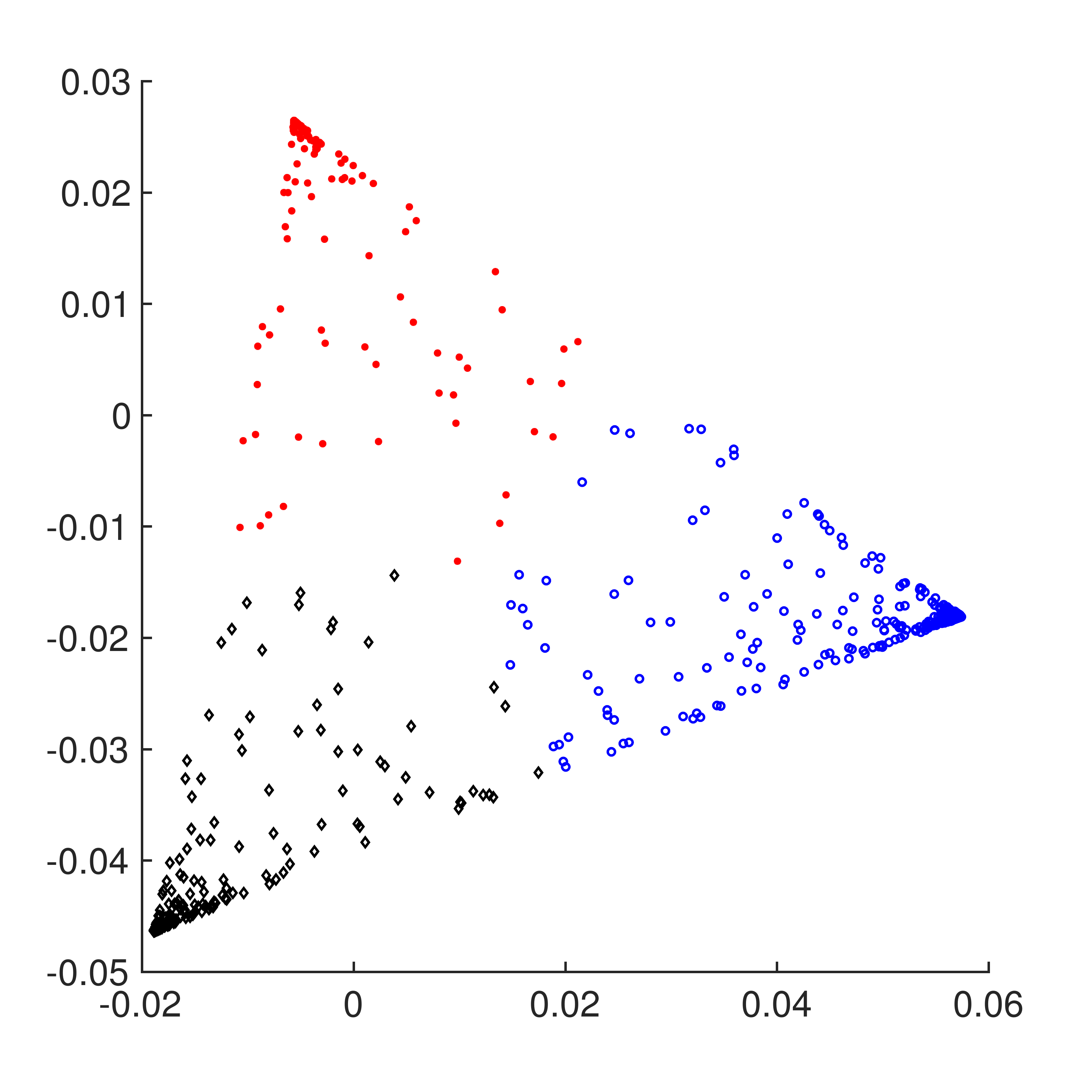}\\
	(c) & (d)
\end{tabular}
\end{center}
\caption{Spectral clustering of synthetic data ($k=3$) with different setup. (a) all synthetic marker genes with low noise (b) all synthetic genes with low noise (c) all synthetic genes with medium noise (d) all synthetic genes with high noise.}%
\label{fig:datatop}%
\end{figure}

Figure \ref{fig:datatop} displays geometric structures of a set of synthetic data. In this case, we take data dimensions $N = 800$, $n=30$, and $k=3$. We define the noise to data ratio (NDR) as $\|\epsilon\|_F/\|{\bf CP}\|_F$ and consider different (low, medium, and high) noise levels in simulations. Eigenvectors of graph Laplacians of these data are computed, and then they are used as coordinates to plot the $N$ points in  Figure \ref{fig:datatop}.   
Note that according to spectral clustering theory, only the second and third eigenvectors are needed (because $k=3$ and the first one is almost a constant vector).
 Figure \ref{fig:datatop} (a) shows the case when $N_1=N_2=300, N_3=200$ while $N_4=0$, i.e., the case of all marker genes. Three different colors represent the three clustered groups. It can be concluded from this subfigure that after clustering, all marker genes will concentrate around the vertices of the $k-1$ simplex since all the corresponding rows of ${\bf C}$ are strongly correlated to some ${\bf e}_l^\top$.  In contrast, Figures \ref{fig:datatop} (b)-(d) show data of all genes, including non-marker genes ($N_1=N_2=N_3=N_4=200$) with low, medium and high levels of noises, respectively. It can be observed that non-marker genes will fill the edge, and strong noises will fill in the interior of the simplex. 
\begin{figure}[h!]
\begin{center}%
\begin{tabular}[c]{ccc}%
	\includegraphics[width=0.33\textwidth]{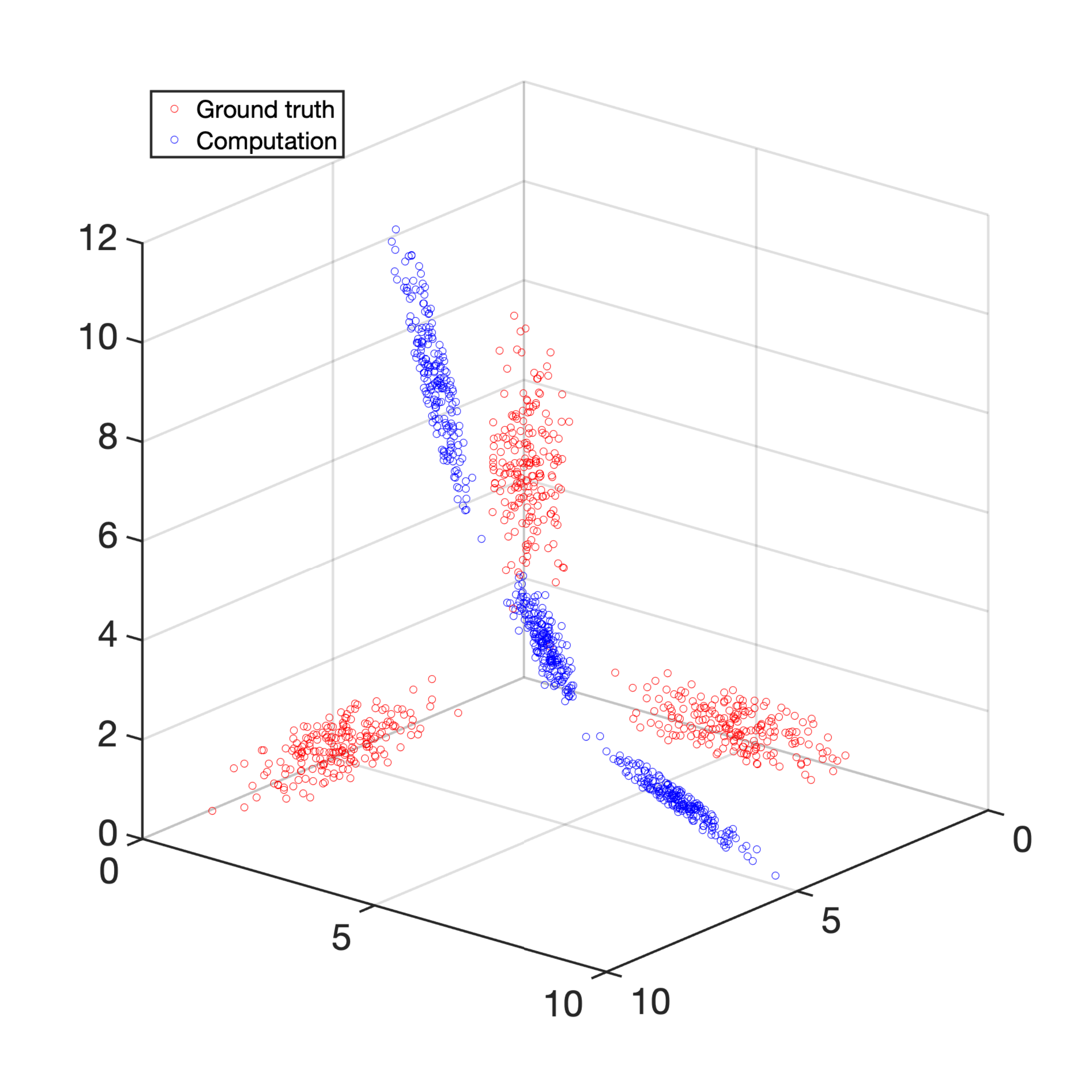} &
	\includegraphics[width=0.33\textwidth]{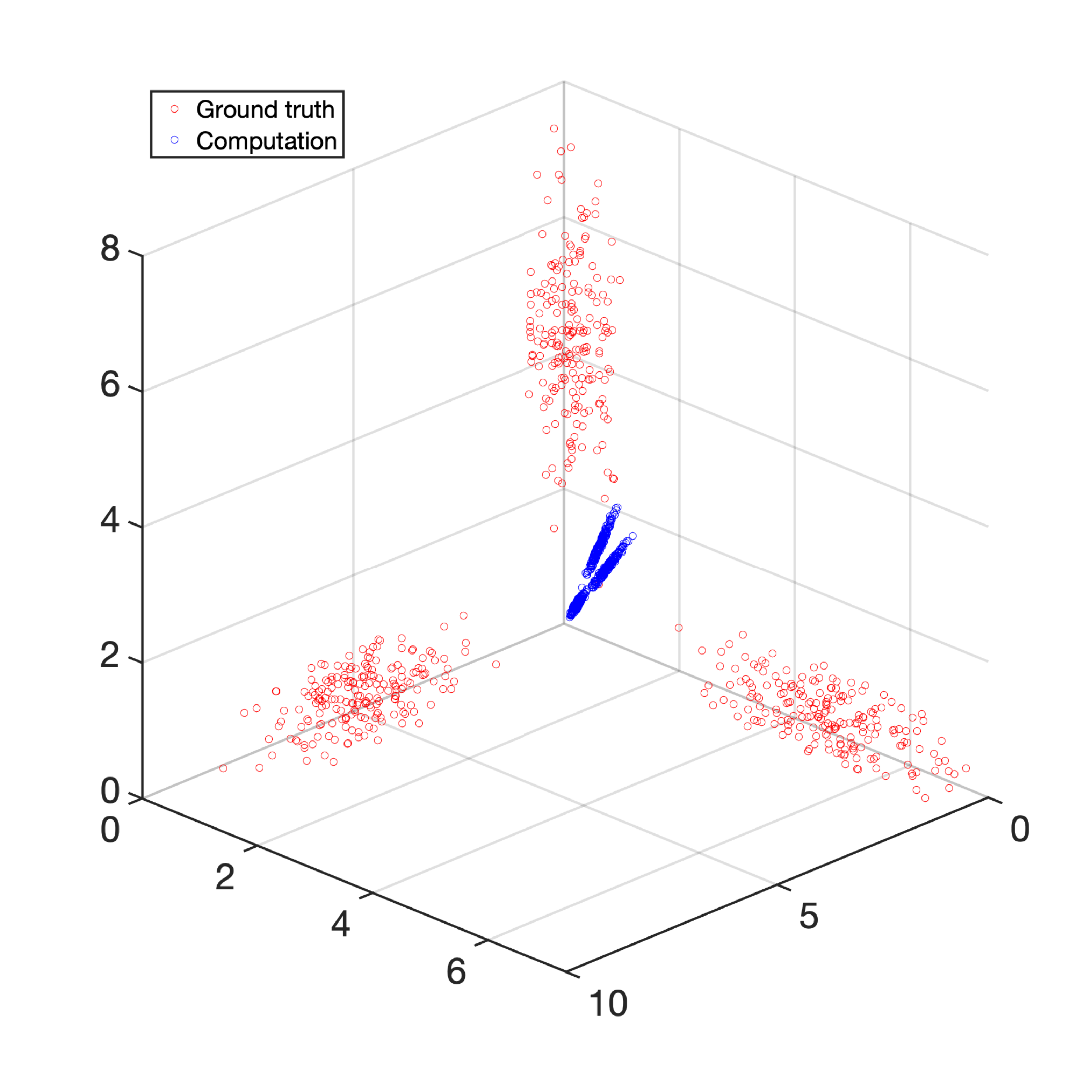} &
	\includegraphics[width=0.33\textwidth]{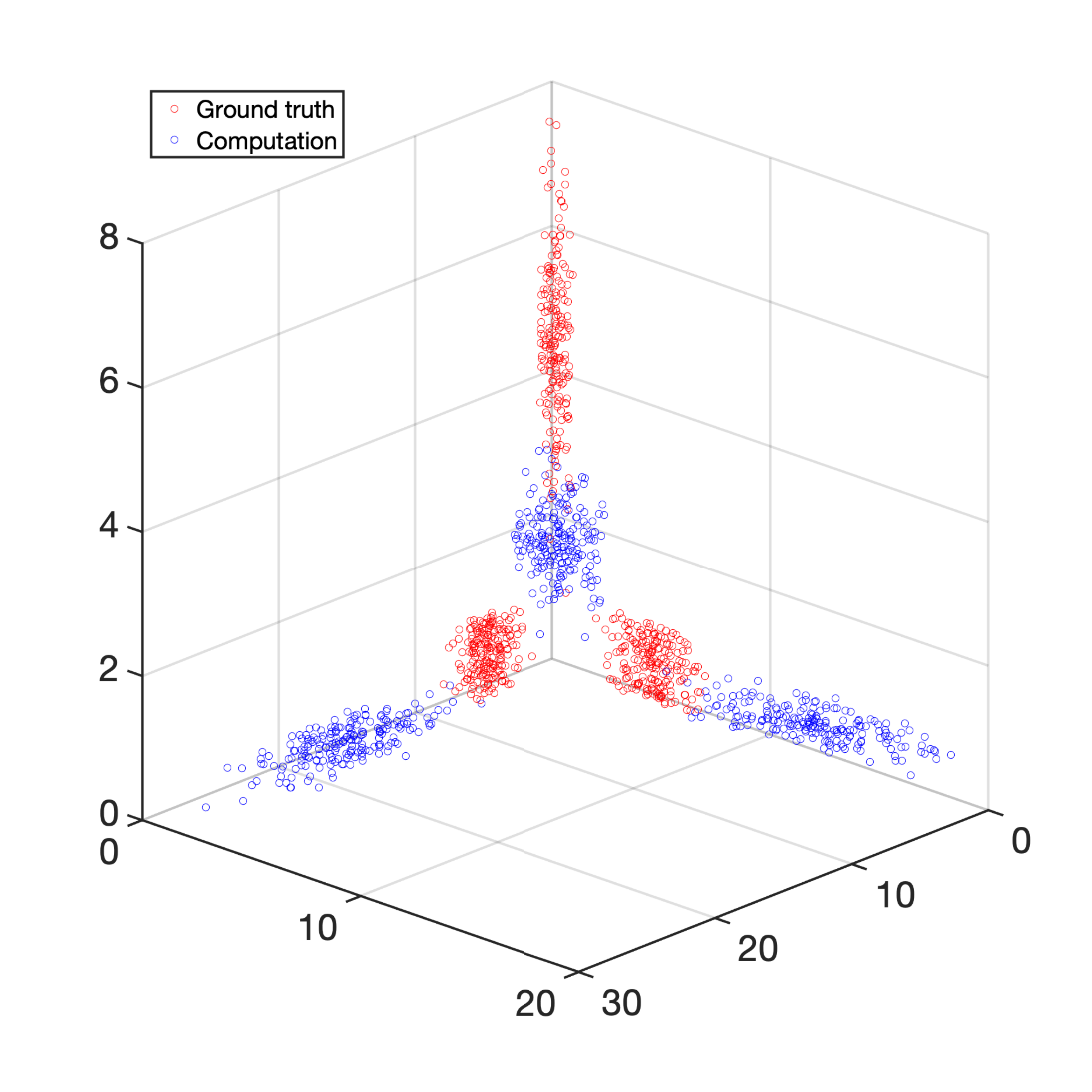}\\
	(i) & (ii) & (iii)
	\end{tabular}
\end{center}
\caption{Complete deconvolution without constraint. (i)-(iii): Comparisons of true and simulated matrices of $\bf C$ with three different initial conditions. }%
\label{fig:noconstraint}%
\end{figure}

In order to show that the proposed constraints are important, we performed the NMF without constraints, by simply setting $\lambda_1=\lambda_2=0$. Initial starting points ${\bf C}_0$ and ${\bf P}_0$ are chosen randomly, so it can be seen in Figure \ref{fig:noconstraint}  that each initial condition will end in different result from others. In these experiments, the stopping criteria are the same ($10^{-5}$) and the relative residues   $\|{\bf G}-{\bf CP}\|_F/\|{\bf G}\|_F$ are the same and consistent to the NDR. Hence we can claim that the approximated stationary points are achieved but
none of them is even close to the ground truth. The different solutions are due to the illposedness of the original NMF problem.
\begin{figure}[h!]
\begin{center}%
\begin{tabular}[c]{cc}%
	\includegraphics[width=0.45\textwidth]{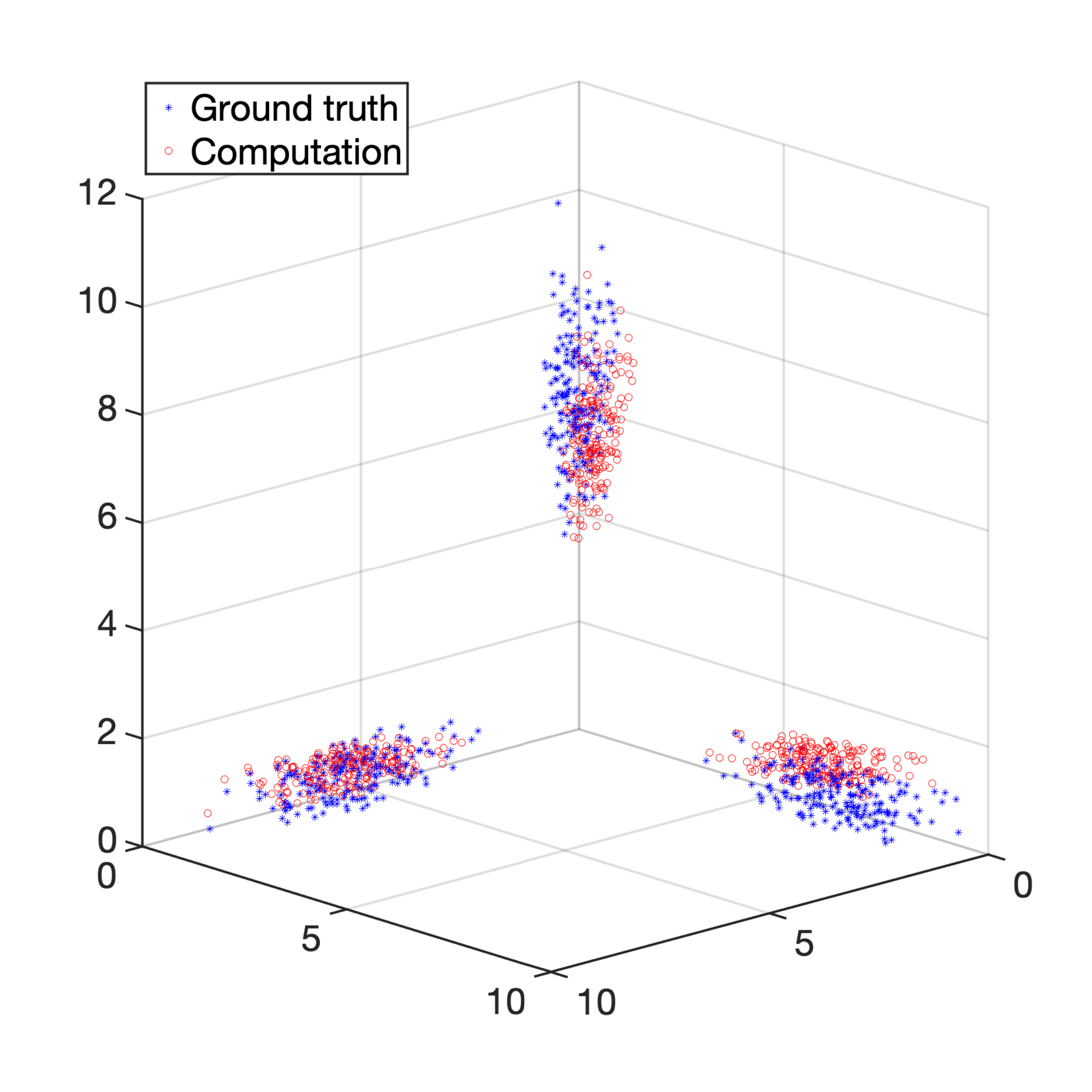} &
	\includegraphics[width=0.45\textwidth]{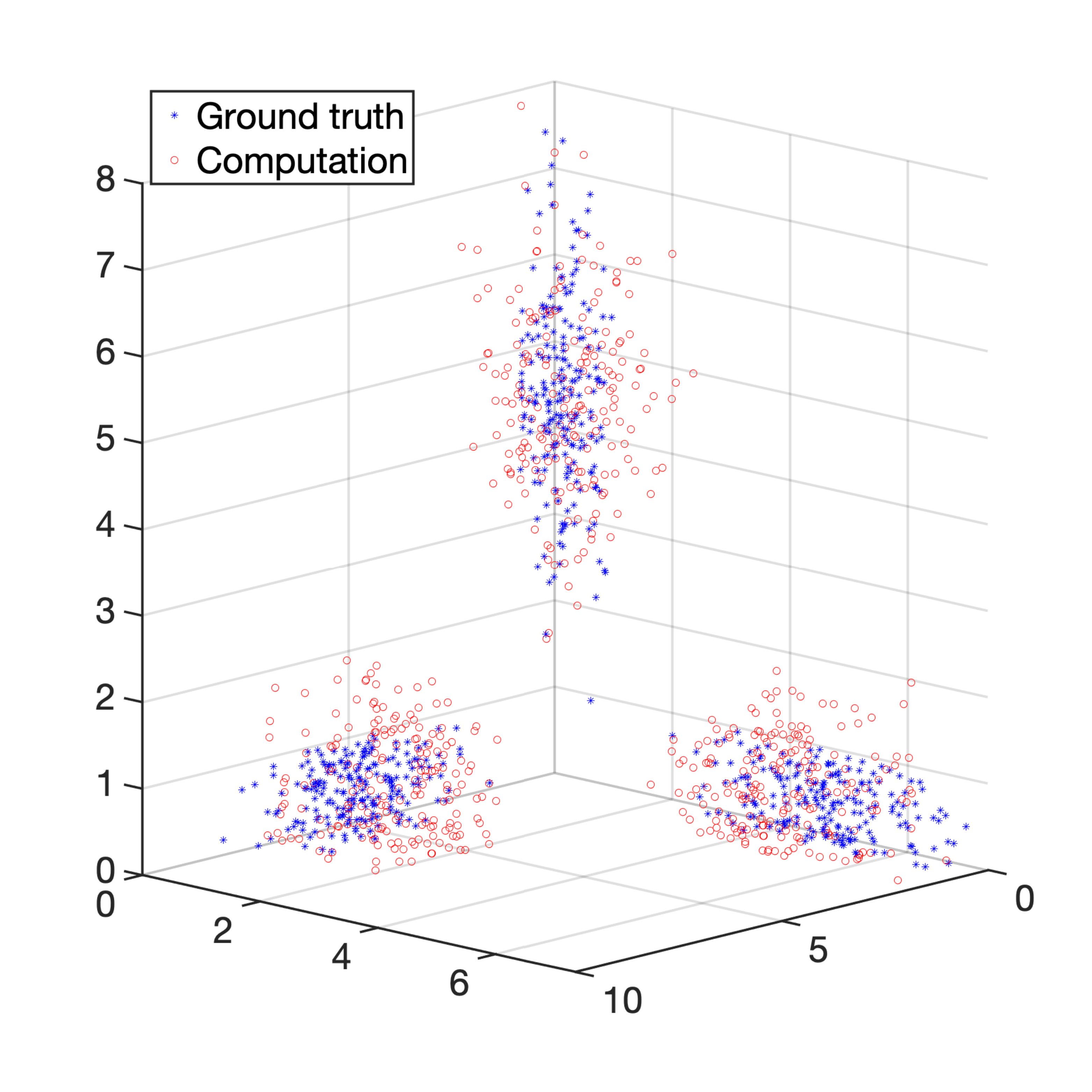}\\
	\includegraphics[width=0.45\textwidth]{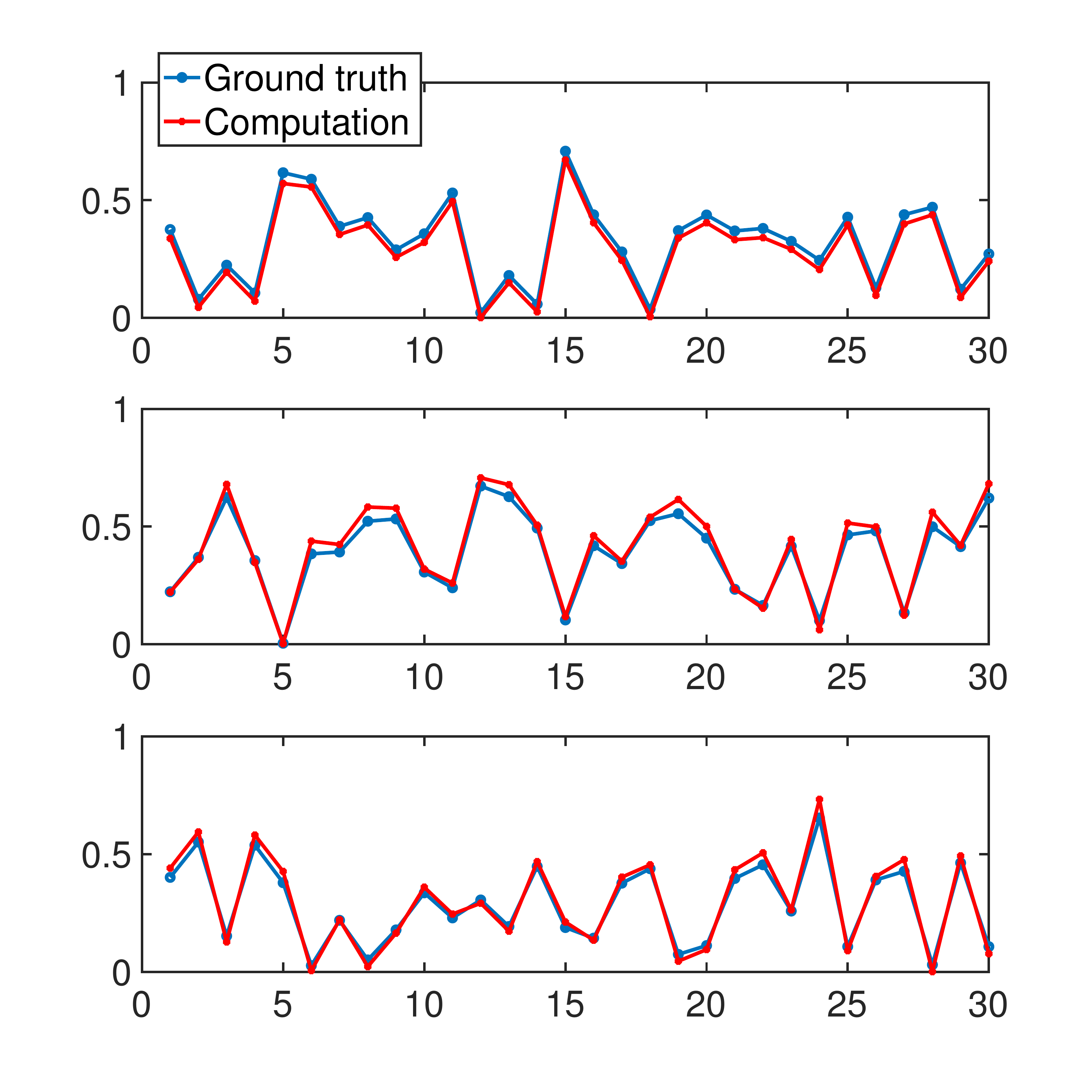} &
	\includegraphics[width=0.45\textwidth]{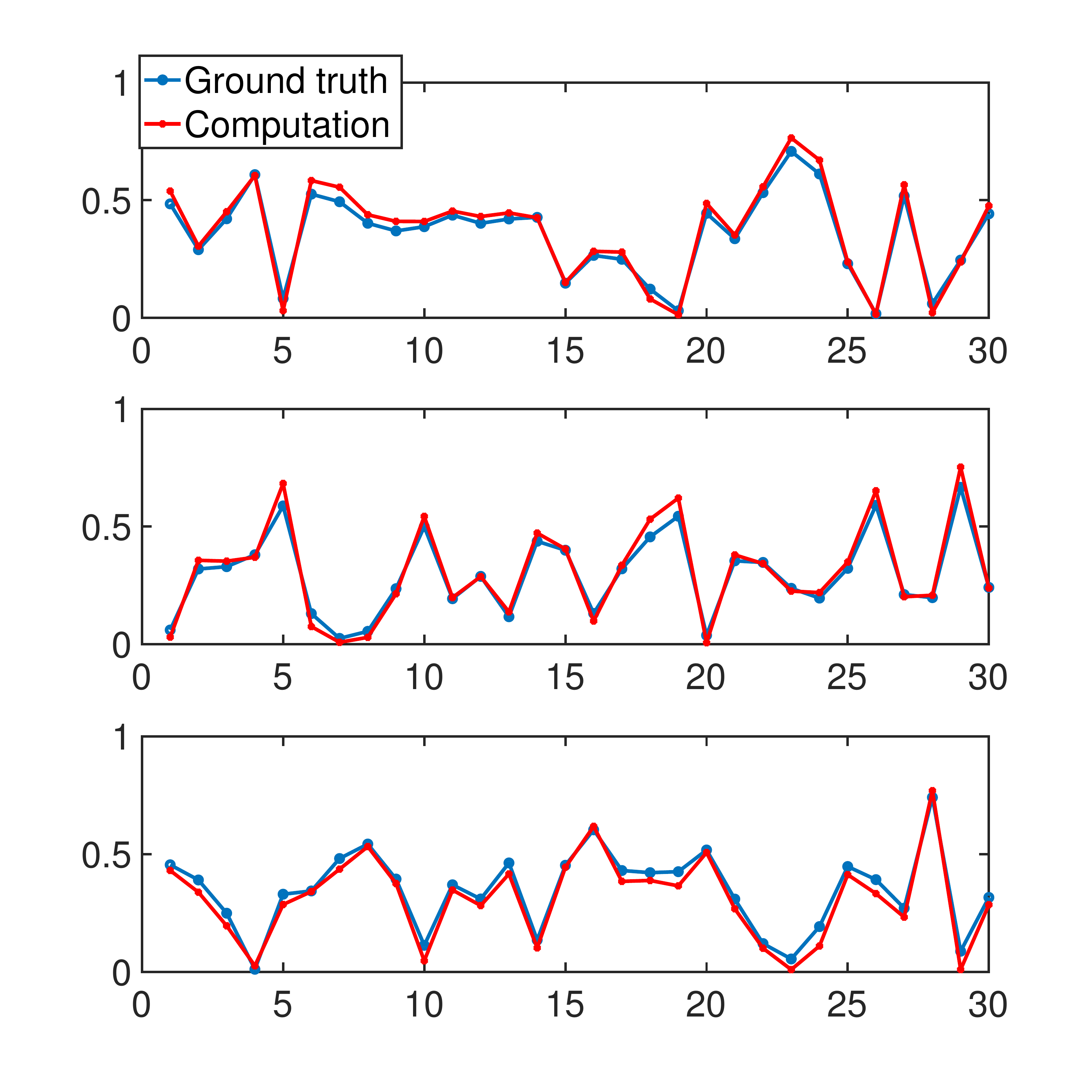}
\end{tabular}
\end{center}
\caption{Comparison of simulations of C and P with ground truth.  Qualitative results are given in Table \ref{table:cp}.}%
\label{fig:CPnoise_nobias}%
\end{figure}

 Figure \ref{fig:CPnoise_nobias} present  computational results of the GS-NMF model on these synthetic data. Two sets  of randomly generated matrices ${\bf C}^*$ and ${\bf P}^*$, with different levels of noises (NDR = $0.071$, $0.336$) are used to obtain data ${\bf G}$. Then comparisons of  ground truth (blue) to the corresponding computational results (red) of the two sets of data are displayed in the left and right panels in Figure \ref{fig:CPnoise_nobias}. The first and second rows are for $\bf C$ and $\bf P$, respectively.  It can be seen that the solutions are remarkably more reasonable from the GS-NMF model. 
\begin{table}[th!]
\centering
\begin{tabular}
[c]{l|c|c|c}\hline\hline
NDR& erros in ${\bf C}$ & errors in ${\bf P}$ & Relative residue \\\hline
0.071 & 0.0901 & 0.0444 &0.0693 \\
0.162 & 0.1007 & 0.0457 &0.1543 \\
0.336 & 0.1372 & 0.0545 &0.3024 \\
0.599 & 0.1667 & 0.0569 &0.4888 \\
\hline\hline
\end{tabular}\label{table:cp}%
\caption{Quantitative results of the GS-NMF with different noise to data ratio (NDR).}%
\end{table}
Quantitative results can be found in Table \ref{table:cp}, where relative errors (comparing to ground truth) of ${\bf C}$, ${\bf P}$, and relative residues $\|{\bf G}-{\bf CP}\|_F/\|{\bf G}\|_F$ are displayed for data with different NDR. As indicated by both Figure \ref{fig:CPnoise_nobias} and the table, errors in matrix ${\bf C}$ increases more obviously when more noises present (larger NDR). On the contrary, computation of matrix ${\bf P}$ seems less vulnerable to noise levels. We observed that the relative residues $\|{\bf G}-{\bf CP}\|_F/\|{\bf G}\|_F$  for all iterations have been already  comparable to the NDR, and this implies that pursing even smaller residues in the cost functions is not necessary.

%=========================results of biased proportion, not present in this work===============================
%\begin{figure}[ptb]
%\begin{center}%
%\begin{tabular}[c]{cc}%
%	\includegraphics[width=0.45\textwidth]{images/Csim_biased1} &
%	\includegraphics[width=0.45\textwidth]{images/Csim_biased2}\\
%	\includegraphics[width=0.45\textwidth]{images/Psim_biased1} &
%	\includegraphics[width=0.45\textwidth]{images/Psim_biased2}
%\end{tabular}
%\end{center}
%\caption{Simulation of C and P with noises. Qualitative results given in a table}%
%\label{fig:CPnoise_nobias}%
%\end{figure}
%=================================================================================================

\begin{figure}
	\begin{center}
		\includegraphics[width=0.75\textwidth]{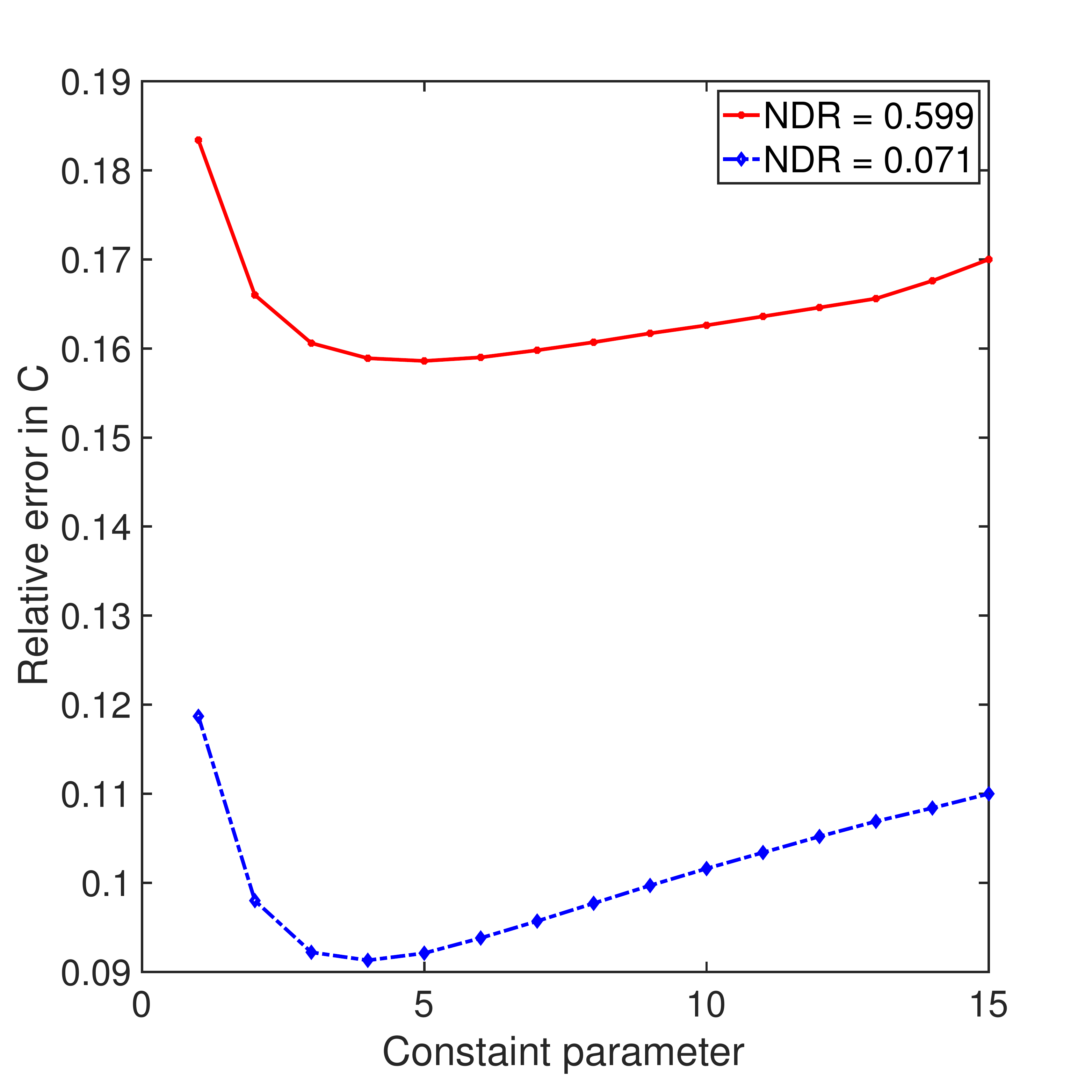} 
	\end{center}
	\caption{Errors in $\bf C$ and parameter $\tilde{\lambda}$}%
\label{fig:pars}%
\end{figure}

The major parameters of the propose algorithm include penalty parameters $\rho$ and $\gamma$, as well as the constraint coefficients $\lambda_1$, $\lambda_2$ in (\ref{eqn:con1}) and (\ref{eqn:con2}). When choosing $\rho$ and $\gamma$,  a grid search is performed with candidates evenly spaced over the interval  $[10^2, 10^5]$. We found that  errors in $\bf C$ and $\bf P$ decrease for larger parameters, while too large values for $\rho$ and $\gamma$ will introduce matrix singularity in the algorithm. Computational results in Table \ref{table:cp} are obtained with $\rho = 1.6\times 10^3$ and $\gamma = 1.5\times 10^4$. For the constraint parameters, we simply take $\lambda_1 = \lambda_2= \lambda$ and rescale to $\tilde{\lambda}=\lambda/\rho$ as it is defined in Eq. (\ref{eqn:scaledADMM}). Empirically, smaller $\tilde{\lambda}$ means less constraints on the geometric structure of ${\bf C}$ hence could damage solution identifiability. On the other side, the extreme case $\tilde{\lambda}\to\infty$ implies the strong identifiability, which is not realistic when noises present. Figure \ref{fig:pars} displays relative errors in $\bf C$, when parameters $\rho$ and $\gamma$ fixed as above but $\tilde{\lambda}$ varies. It can be concluded that for both noise levels (NDR = 0.599 and 0.071), the change of errors against $\tilde{\lambda}$ is not monotone. For the testing data, $\tilde{\lambda}\approx 4$ or $5$ seems the best choice for computational accuracy. How to chose reasonable parameter $\tilde{\lambda}$ according to different data set could be a future study.

\subsection{Algorithm results on biological data}

We also validate the proposed algorithms by realistic biological data from GSE19830 \cite{shen2010cell}. This data set is obtained from tissue samples of the brain, liver, and lung of  a single rat using expression arrays (Affymetrix). Homogenates of these three type of tissues were mixed together at the cRNA homogenate level with a known proportion, and then the gene expression pattern of every mixed sample was measured. The GSE19830 data set mimics the common scenario of heterogeneous biological samples which vary in the relative frequency of the component subsets from one to another and has been used in a few literatures \cite{kang2019cdseq,zaitsev2019complete} to validate computational algorithms. For this dataset, we know cell type $k=3$ and tissue sample number $n=33$. After necessary data preprocessing to exclude obvious outliners (row norm, column norm, etc), we take $N=10,000$ out of $\sim 12,000$ total genes. 

\begin{figure}[ptb]
\begin{center}%
\begin{tabular}[c]{cc}%
	\includegraphics[width=0.5\textwidth]{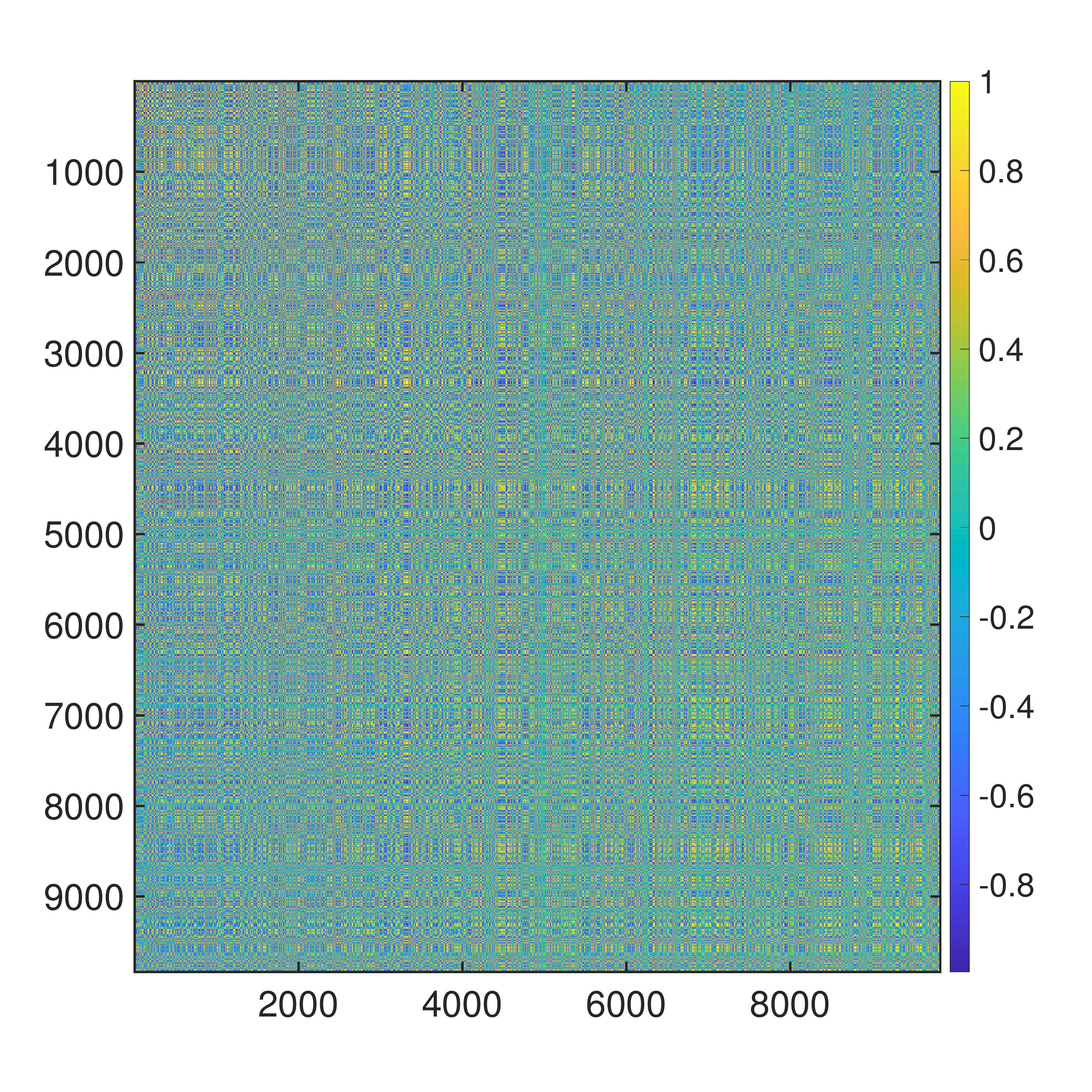} &
	\includegraphics[width=0.5\textwidth]{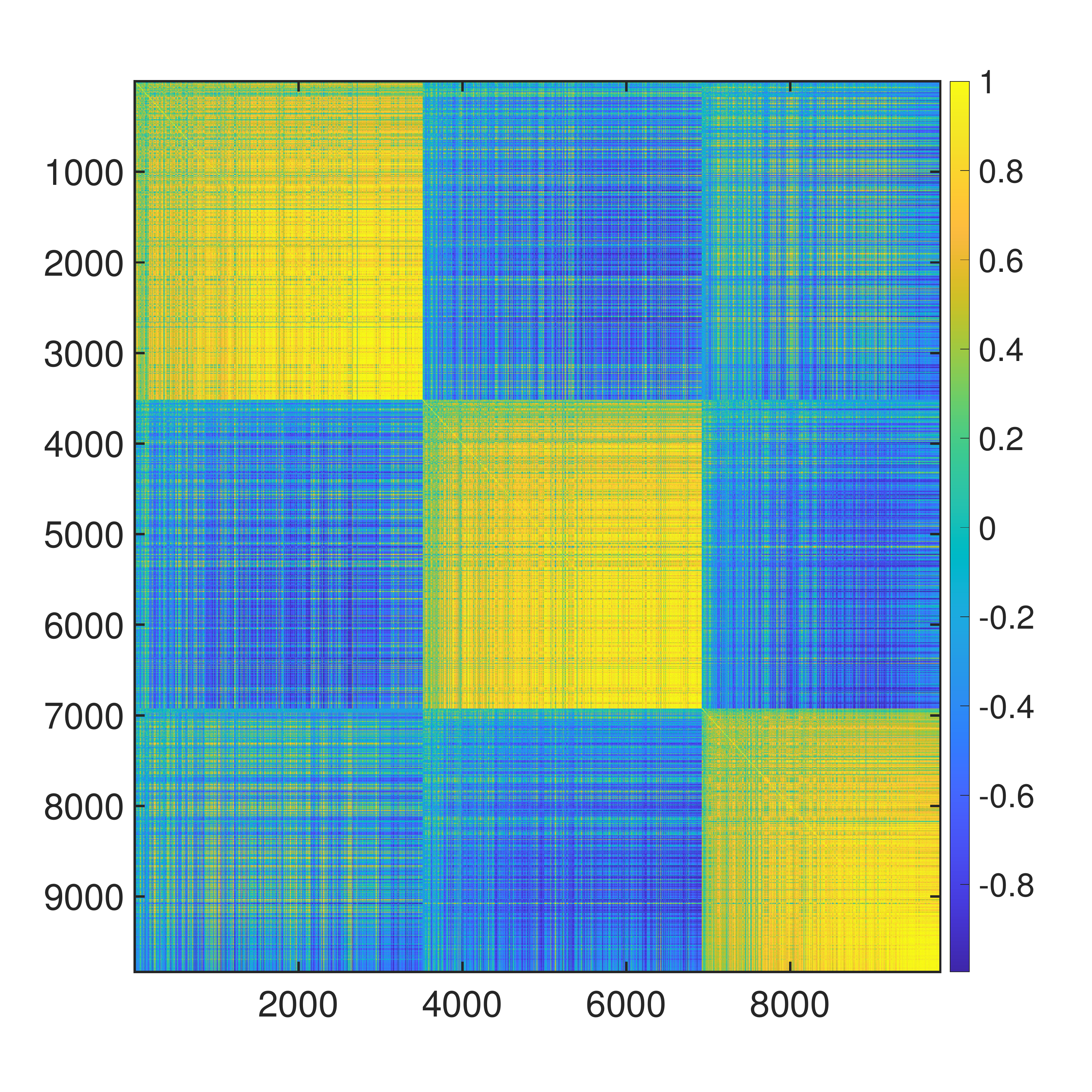}\\
	unclustered & clustered
	\end{tabular}
\end{center}
\caption{Correlations of gene expressions across tissue samples.}%
\label{fig:corr}%
\end{figure}

Figure \ref{fig:corr} displays the structure of data ${\bf G}\in \R^{N\times n}$: mutual correlations of the rows of ${\bf G}$, i.e. gene expression of the $10,000$ genes in those 33 samples, are computed and shown as  heat maps, before (left) and after (right) clustering/permutation. Data $\bf G$ has a clear structure: there are three clusters and  within each clustered group, gene expressions are strongly correlated. Based this clustering and evaluation of correlations, marker genes could be computationally identified.

\begin{figure}[ptb]
\begin{center}%
\begin{tabular}[c]{cc}%
	\includegraphics[width=0.45\textwidth]{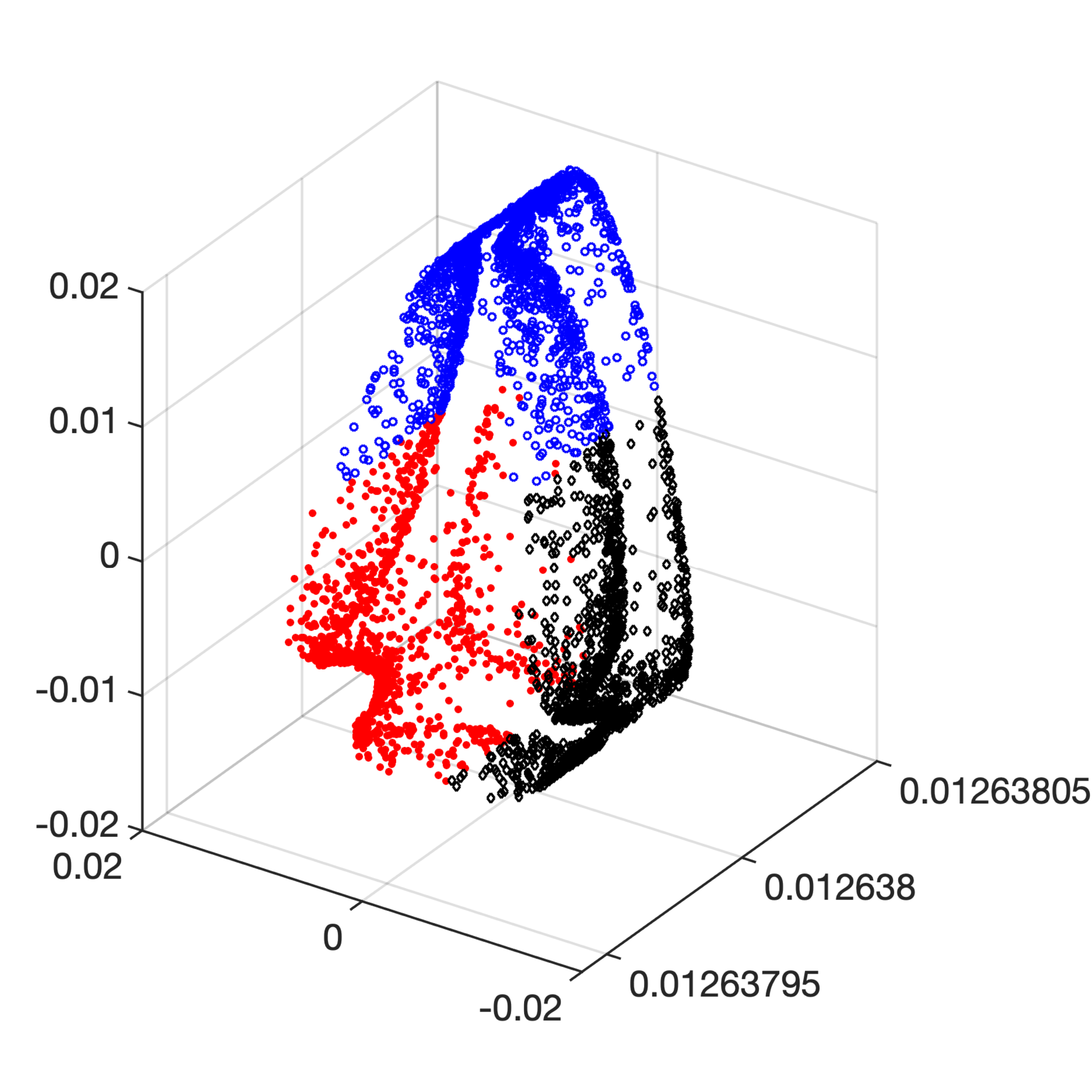} &
	\includegraphics[width=0.45\textwidth]{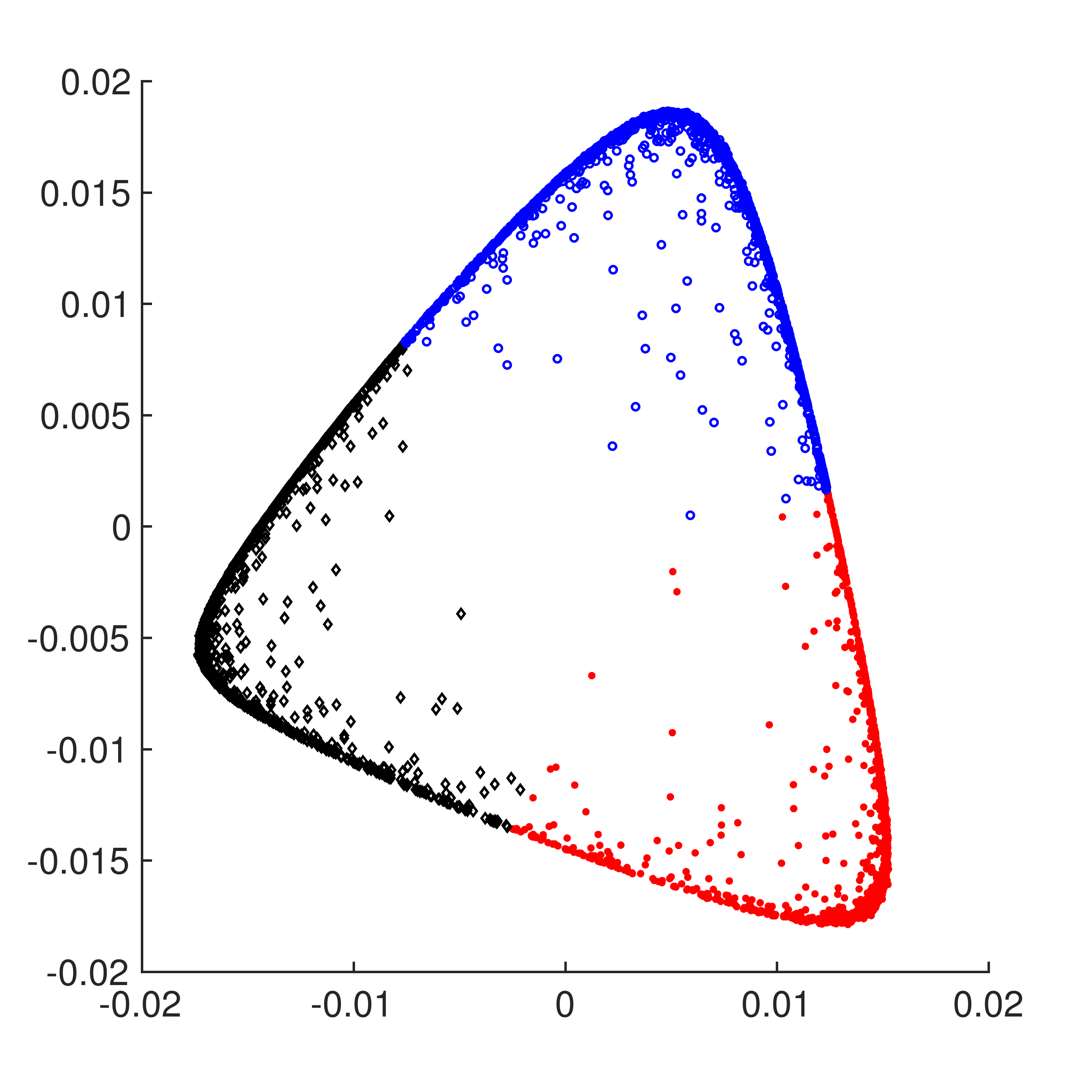}
\end{tabular}
\end{center}
\caption{3D (left) and 2D (right) presentation of spectral clustering of GSE19830 data.}%
\label{fig:3celldata_2D3D}%
\end{figure}

\begin{figure}[ptb]
\begin{center}%
\begin{tabular}[c]{cccc}%
	\includegraphics[width=0.23\textwidth]{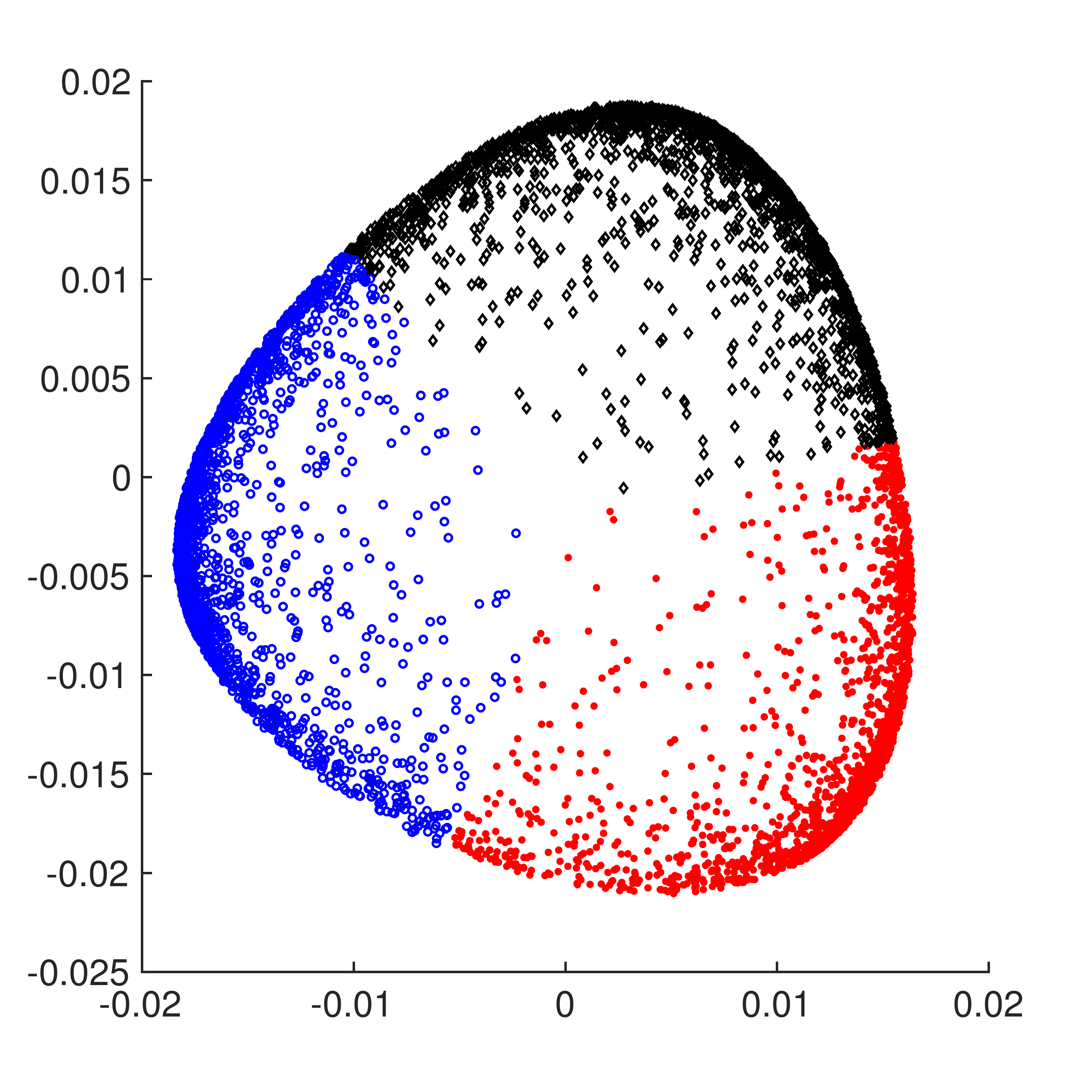}&
	\includegraphics[width=0.23\textwidth]{images/3cells2D_5}&
	\includegraphics[width=0.23\textwidth]{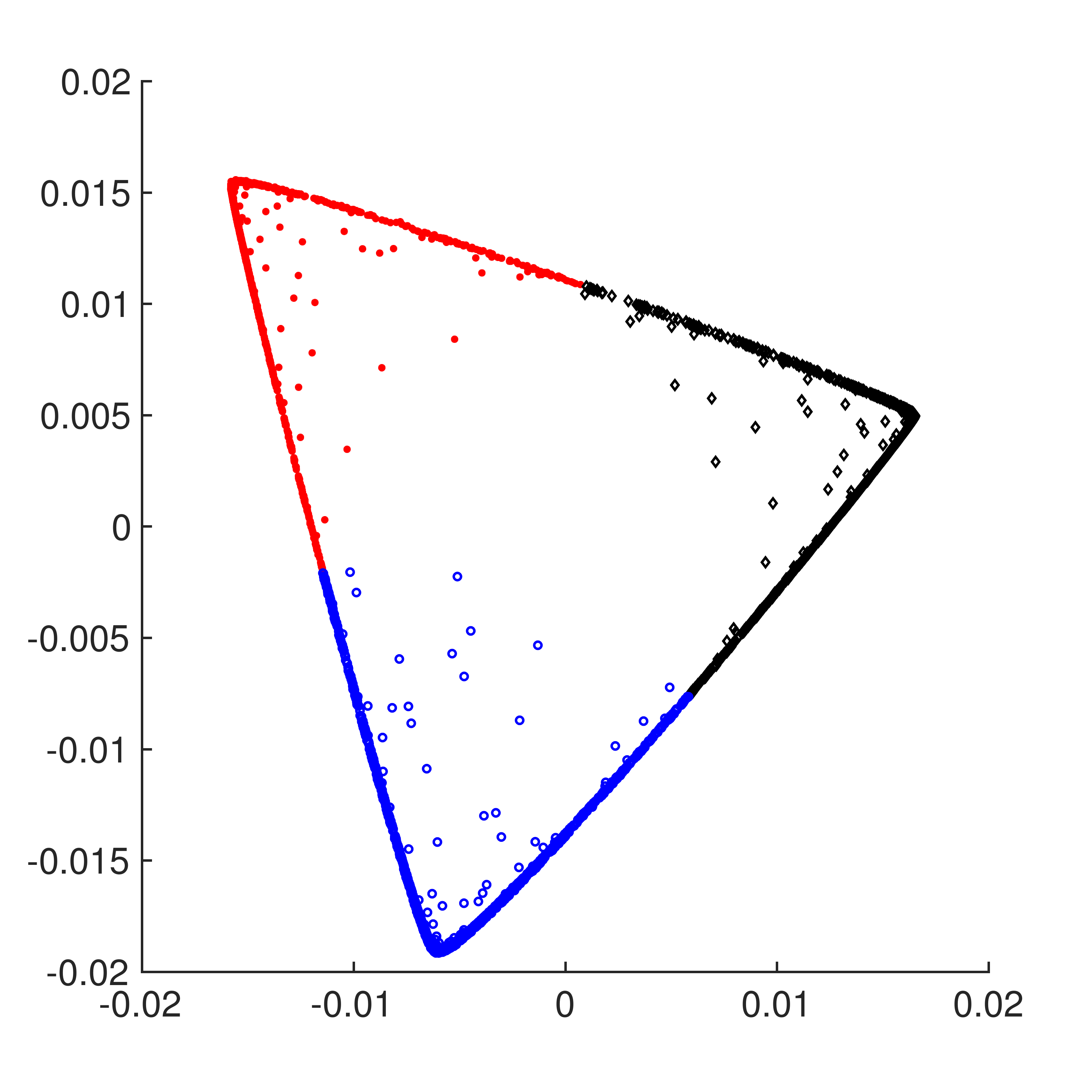}&
	\includegraphics[width=0.23\textwidth]{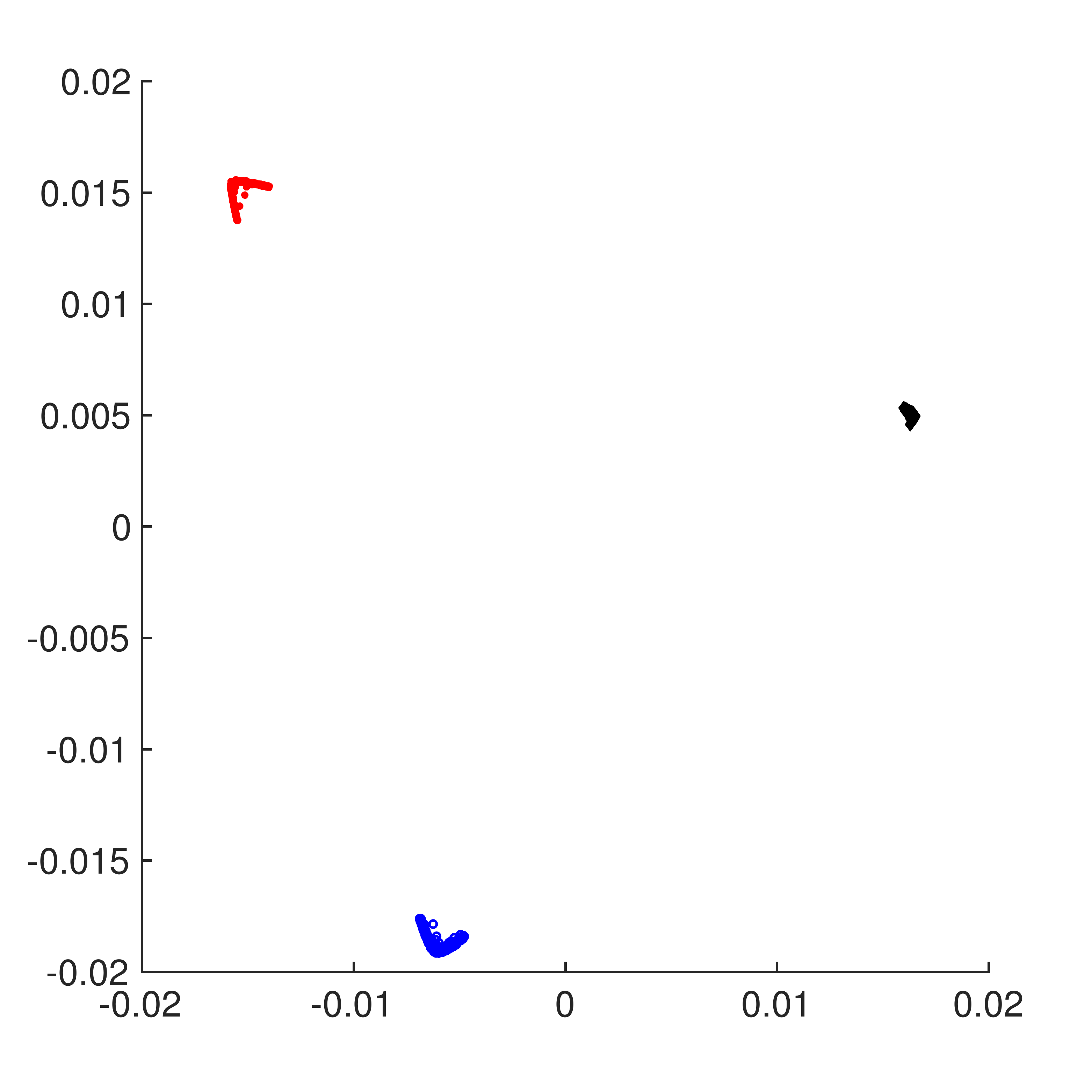}\\
	(a) $\sigma = 1.0$ & (b) $\sigma = 0.5$ & (c) $\sigma= 0.2$ & (d) $\sigma = 0.2$,  $\mathcal{S}_r$
\end{tabular}
\end{center}
\caption{Spectral clustering of GSE19830 data with $\sigma = 0.2, 0.5, 1$. (a)-(c): whole set; (d) Set $\mathcal{S}_r$}%
\label{fig:3celldata_2D}%
\end{figure}

Figures \ref{fig:3celldata_2D3D}-\ref{fig:3celldata_2D} show some clustering details. 
Since the clustering is based on eigenvectors of the graph Laplacian of data ${\bf G}$, we plot the first three eigenvectors of ${\bf L}$, with each column as the $x$-, $y$-, and $z$-coordinates of the $N$ dots, respectively  in the left panel of  Figure \ref{fig:3celldata_2D3D}. Different colors represent the clustered groups. The first eigenvector is almost a constant vector, so it is convenient to just display the 2D data for the rest of figures, as in the right panel. All the dots are distributed in a $(k-1)$-simplex, and each cluster is identified with its vertex. 

Notice that exploring the data structure  depends on the parameter $\sigma$ in Eq. (\ref{eqn:omga}). Figure \ref{fig:3celldata_2D} shows data distribution with values of $\sigma = 1, 0.5$, and $0.2$. From Eq. (\ref{eqn:omga}) we see that the connectivity of any two vertices in the graph increases for large value of $\sigma$. This feature is displayed in Figure \ref{fig:3celldata_2D} (a), (b) and (c): when $\sigma = 1$, all data points are distributed within almost a circle and the clustering is not that significant. On  the other hand, when $\sigma = 0.2$, the three vertices of the triangles naturally define the three clusters. In our experiment, the graph loses majority of connectivity for even smaller value, so $\sigma = 0.2$ is used for all the simulations.  To determine marker genes, we pick a subset $\mathcal{S}_r$ from each colored group and they are chosen as the one that are closest to each vertex. With $\sigma = 0.2$, each set  $\mathcal{S}_r, r = 1, 2, 3$ contains 1,000 entries, and the corresponding data points are shown in Figure \ref{fig:3celldata_2D} (d).

\begin{figure}[ht!]
	\begin{center}
		\includegraphics[width=0.95\textwidth]{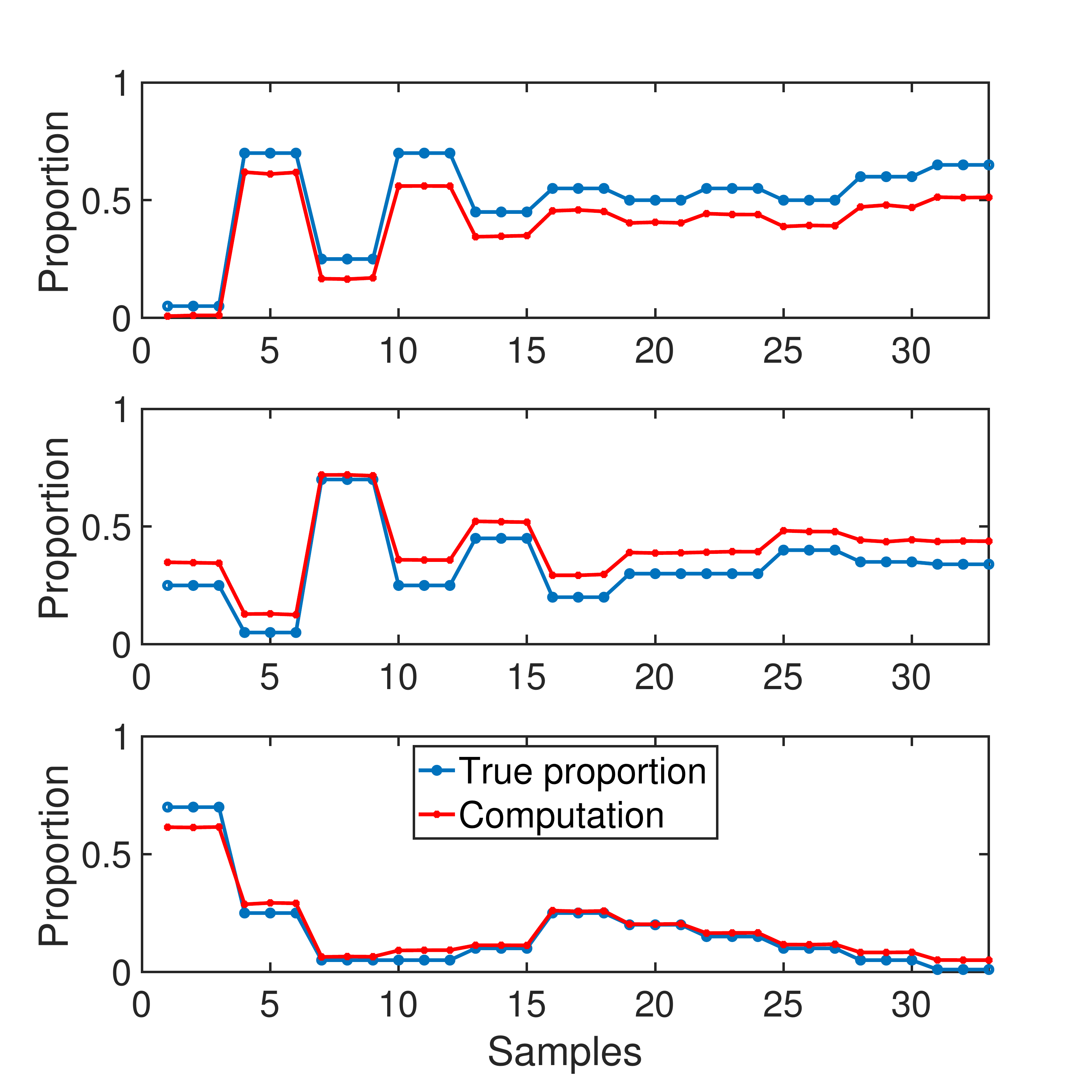} 
	\end{center}
	\caption{Comparison of simulated and true cellular proportions. From top to bottom are for liver, brain, and lung cells.}%
\label{fig:cellP}%
\end{figure}

With such parameter choices and set $\tilde{\lambda}=0.6$, we applied the proposed algorithms to data set GSE19830. Figure \ref{fig:cellP} shows the comparison between computed (red) cellular composition (liver, brain, lung, from top to bottom) in bulk tissue samples and ground truth (blue). In  the 33 samples, 11 different cellular compositions were used and each of them was replicated three times. The computational results have reproduced this pattern.  Additionally, the simulated cellular proportions fit the ground truth fairly well, especially for the third cell type. Note that correlations of the blue/red curves in the three panels of Figure \ref{fig:cellP} are $0.9916, 0.9916$ and $0.9997$. It indicates that simulation $\tilde{\bf P}$ and  ground truth $\bf P$  differ by merely a scaling factor, i.e. $\tilde{\bf P}=\text{diag}(s_1, s_2,s_3){\bf P}$. This phenomenon is majorly due to the definition of uniqueness of the NMF in Eq. (\ref{eqn:def}). A future direction could be  improvement of the algorithm, in order to make the scaling diagonal matrix close to an identity matrix. 
%\begin{figure}[ptb]
%\begin{center}%
%\begin{tabular}[c]{cc}%
%	\includegraphics[width=0.5\textwidth]{images/3cells2D_2} &
%	\includegraphics[width=0.5\textwidth]{images/3cells2D_2_m}
%	\end{tabular}
%\end{center}
%\caption{Spectral clustering of GSE19830 data with $\sigma = 0.2$ (left)and subset $\mathcal{S}_r$ (right).}%
%\label{fig:marker}%
%\end{figure}

\section{Conclusions} \label{sec:con}

With current technologies, large-scale bulk-RNAseq data are available to study molecular pathways implicated in various diseases. Differences of transcriptome-wide gene expression profiles (GEPs) among patients and controls will reveal novel insights into genes and pathways, so they are potentially helpful for drug targets therapeutics.  However, it is always a challenge whether disease-associated GEP differences in tissues are due to changes in cellular composition of tissue samples, or due to GEP changes in specific cells. Although single-cell RNAseq data can be used or serve as references, such approaches remain costly, cumbersome and limited in sample size. In contrast, computational approaches can be used to decompose the more reliable bulk-tissue RNAseq data.
  	 In this paper we  develop a robust mathematical model and corresponding  computational  algorithms for complete data deconvolution.  The major technique is nonnegative matrix factorization (NMF), which has a wide-range of applications  in the machine learning community. Meanwhile, the NMF is a well-known strongly ill-posed problem, so a direct application of it to RNAseq data will suffer severe difficulties in the interpretability of solutions. To address this issue, we leverage the biological concept of marker genes, combine it with the solvability conditions of the NMF theories, and hence  develop a geometric structured guided optimization problem. In this approach,  the geometric structure of bulk tissue data is first explored by the spectral clustering technique. In this step, correlations graph among GEPs across tissue samples is established, and more importantly, marker genes for each cell types are identified.   Then,  information of marker genes  is integrated as solvability constraints, while the overall correlation graph is used as manifold regularization. The resulting non-convex optimization problem, termed as geometric structured nonnegative matrix factorization (GS-NMF) model is numerically solved under the framework of  alternating direction method of multipliers (ADMM). Finally, synthetic and biological data are used to validate the proposed model and algorithms. With this novel method, solution interpretability is significantly improved  and accuracy is satisfactory comparing to the ground truth for both types of data. It is worthwhile to note that all simulation results may still suffer a linear scaling factor comparing to the ground truth. Unfortunately this is nothing to do with marker gene selection, parameter choices, or algorithm accuracy,  but due to the inherent definition of NMF solution uniqueness. In the future research, we will combine necessary biological information in realistic applications, to reduce this  scaling ambiguity as much as possible.

\bibliographystyle{ieeetr}
\bibliography{refs}

\end{document}